\documentclass[11pt]{amsart}

\usepackage[utf8]{inputenc}
\usepackage[T1]{fontenc}
\usepackage[english]{babel}

\usepackage{amsmath, amssymb, amsfonts, amsthm, amsbsy}
\usepackage{dsfont}

\usepackage{graphicx}
\usepackage{tikz}
\usepackage{caption}
\usepackage{subcaption}
\usepackage{float}

\usepackage{enumitem}
\usepackage{relsize}
\usepackage{stackengine}

\usepackage{hyperref}
\usepackage[capitalise,compress]{cleveref}
\usepackage{cite}

\usepackage{color}

\usepackage{comment}
\usepackage{psfrag}

\usepackage{etoolbox}

\numberwithin{equation}{section}

\newtheorem{theorem}{Theorem}[section]
\newtheorem{lemma}[theorem]{Lemma}
\newtheorem{proposition}[theorem]{Proposition}
\newtheorem{corollary}[theorem]{Corollary}

\theoremstyle{definition}
\newtheorem{definition}[theorem]{Definition}

\newtheorem{con}[theorem]{Conditions}

\theoremstyle{remark}
\newtheorem{remark}[theorem]{Remark}

\newcommand{\RR}{\mathbb{R}}

\newcommand{\BB}{\mathbb{B}}
\newcommand{\HH}{\mathbb{H}}
\newcommand{\SPH}{\mathbb{S}}
\newcommand{\RP}{\mathbb{RP}}

\newcommand{\proj}{\operatorname{proj}}
\newcommand{\Gnom}{\operatorname{Gnom}}

\newcommand{\conv}{\operatorname{conv}}
\newcommand{\pow}{\operatorname{pow}}

\newcommand{\Vol}{\operatorname{Vol}}
\newcommand{\dd}{\mathrm{d}}

\newcommand{\ctime}{\mathcal{C}}
\newcommand{\loren}[2]{\langle #1, #2 \rangle_{L}}

\newcommand{\xx}{\tilde{x}}
\newcommand{\yy}{\tilde{y}}

\newcommand\barbelow[1]{\stackunder[1.2pt]{$#1$}{\rule{.8ex}{.075ex}}}

\DeclareMathOperator{\arcosh}{arcosh}

\date{\today}

\makeatletter

\renewcommand{\l@section}{%
  \@tocline{1}{6pt}{0pt}{2.5em}{\bfseries}}

\renewcommand{\l@subsection}{%
  \@tocline{2}{0pt}{1.5em}{3em}{}}

\patchcmd{\@tocline}
  {\hfil}
  {\ifnum#1=2\relax\dotfill\else\hfil\fi}
  {}
  {\PackageError{toc-format}{Could not patch amsart contents}{}}

\makeatother

\begin{document}

\sloppy

\title[Hyperbolic Laguerre tessellations]{Random Laguerre tessellations in hyperbolic space} 

\author{Anna Gusakova}
\address{Anna Gusakova, Universität Münster, Institut für Mathematische Stochastik, Münster, Germany}
\email{gusakova@uni-muenster.de}

 \author{Panagiotis Spanos}
\address{Panagiotis Spanos, Dipartimento di Matematica,
Sapienza Universit\`a di Roma, Roma, Italy}
\email{panagiotis.spanos@uniroma1.it}

\date{\today}

\renewcommand\footnotemark{}

\begin{abstract}
We introduce Laguerre tessellations in hyperbolic space as a weighted generalization of Voronoi tessellations and a hyperbolic analogue of Euclidean power diagrams. We present geometric and analytic constructions and establish conditions under which a collection of weighted points generates a well-defined tessellation. For Poisson point processes of weighted points, we derive necessary and sufficient conditions on the intensity measure ensuring that the resulting diagram is almost surely a tessellation. We then prove that hyperbolic Poisson-Laguerre tessellations are stable under intersections with totally geodesic subspaces, with an explicit transformation of the intensity measure. As applications, we identify the Poisson-Laguerre tessellations arising as sections of hyperbolic Poisson-Voronoi and ideal Poisson-Voronoi tessellations. In particular, the latter have unbounded cells, while their proper totally geodesic sections have bounded cells almost surely.
\end{abstract}
\setcounter{tocdepth}{2}

\maketitle


\baselineskip=16pt

\section{Introduction and main results}

A Voronoi tessellation of a metric space is generated by a collection of points,
where each cell consists of all points, which are at least as close to its generator as to any other generating point. When the generators form a homogeneous Poisson point process of
intensity $\lambda>0$, the resulting Poisson-Voronoi tessellation is one of
the fundamental models of stochastic geometry and has been extensively studied
in Euclidean space (see, e.g., \cite{Stoyan,Moller,SchneiderWeil}).

Random tessellations beyond the Euclidean setting have recently attracted
considerable attention, in particular in hyperbolic spaces and on Riemannian
manifolds \cite{Isokawa2000,BenjaminiPaquettePfeffer2018,
CalkaChapronEnriquez2021,HeroldHugThaele2021}. The exponential volume growth
and lack of scale invariance of hyperbolic space $\mathbb{H}^d$ give rise to
phenomena, which have no Euclidean counterparts. A striking example is the existence of nontrivial low-intensity limit of hyperbolic Poisson-Voronoi tessellation. Namely, as the intensity $\lambda$ of underlying Poisson point process tends to zero, its Voronoi tessellation converges in
distribution to a nontrivial isometry-invariant tessellation consisting of
countably many unbounded cells, which is called \emph{ideal Poisson-Voronoi tessellation} \cite{DAchilleEtAl2026}. The ideal Poisson-Voronoi
tessellation has been used to obtain bounds for Cheeger constants of hyperbolic surfaces \cite{BudzinskiCurienPetri2025}, while related constructions in
higher-rank symmetric spaces have been used to prove the fixed-price-one
property for lattices in higher-rank semisimple Lie groups
\cite{FraczykMellickWilkens2025}. Further recent work has determined the face
volume densities of hyperbolic Poisson-Voronoi and ideal Poisson-Voronoi tessellations
\cite{DAchilleThaele2026}. These developments suggest that random
tessellations encode important geometric properties of their ambient spaces.
The aim of this paper is to develop a class of random \emph{hyperbolic
Laguerre tessellations} that includes ordinary, ideal and sectional hyperbolic
Poisson-Voronoi tessellations and provides a common framework for studying
their properties.

\subsection{Laguerre tessellations in Euclidean space}

A Laguerre tessellation, also called a \emph{power diagram}, generalizes a
Voronoi tessellation by assigning a weight to each generator
\cite{Aurenhammer1987,Okabe}. For a weighted point $(v,h)\in\RR^d\times\RR$
and $u\in\RR^d$, define the Euclidean power function by
\begin{equation}\label{eq:power}
  \pow_{\rm Euc}((v,h),u):=\|u-v\|^2+h.
\end{equation}
For a countable set $A\subset\RR^d\times\RR$, the cell associated with $(v,h)\in A$ is
\[
  C_{\rm Euc}((v,h),A)
  :=\bigl\{u\in\RR^d:\ \pow_{\rm Euc}((v,h),u)
   \leq\pow_{\rm Euc}((v',h'),u)
   \text{ for all }(v',h')\in A\bigr\}.
\]
Unlike ordinary Voronoi cells, Laguerre cells may be empty. They are convex,
since differences of power functions are affine, and, in particular, for finite generating set they are convex polyhedra. The Laguerre diagram is the family of cells with nonempty interior. The power function can also be interpreted as the arrival
time (at $u$) of a crystal that starts growing at time $h$ from $v$, with radius
$\sqrt{t-h}$ at time $t\geq h$. A complementary geometric description uses
lifting and polarity to relate Laguerre cells to faces of locally polyhedral set
in $\RR^{d+1}$ (see \cite[Chapter~17]{Boisso}).

Random Laguerre tessellations generated by independently marked homogeneous Poisson
process were systematically studied in \cite{LautensackZuyev2008}.
Poisson-Laguerre tessellations, generated by a Poisson point process with intensity measure
$f(h)\,\mathrm dv\,\mathrm dh$ on $\RR^d\times\RR$ were investigated in
\cite{GusakovaInWoldeLuebke2025}. Important examples include the
$\beta$-, $\beta'$- and Gaussian-Voronoi models, with height densities
proportional to $h^\beta\mathbf1_{\{h>0\}}$,
$(-h)^{-\beta}\mathbf1_{\{h<0\}}$ and $e^{\gamma h}$, respectively, in
their admissible parameter ranges
\cite{GusakovaKabluchkoThaele2022,GusakovaKabluchkoThaeleGaussian2022}. The $\beta$-Voronoi model appears naturally when studying sectional Poisson-Voronoi tessellation, namely 
the intersection of a $d$-dimensional Euclidean Poisson-Voronoi
tessellation with a $k$-dimensional affine subspace is a
$\beta$-Voronoi tessellation with $\beta=(d-k)/2-1$
\cite{GusakovaKabluchkoThaele2024}. More generally, sections of
Poisson-Laguerre tessellations are again Poisson-Laguerre tessellations,
with the height density transformed by a Riemann-Liouville fractional
integral of order $(d-k)/2$ \cite{GusakovaInWoldeLuebke2025}.
Related limit theorems have recently been established in
\cite{GusakovaInWoldeLuebkeLimits2026}.
Another link with hyperbolic geometry is that the $\beta'$-Voronoi model
with parameter $\beta=d$ in $\RR^{d-1}$ describes a projection of the
boundary of the typical ideal Poisson-Voronoi cell in upper halfspace model of $\mathbb{H}^d$
\cite{DAchilleEtAl2026}.

Laguerre cells also arise in semi-discrete optimal transport for quadratic
cost, where their weights encode prescribed cell masses
\cite{AurenhammerHoffmannAronov1998,Merigot2011}. Related stochastic optimal
transport constructions and extensions to Riemannian manifolds provide
further motivation for studying weighted tessellations
\cite{HuesmannSturm2013,Huesmann2016}.

\subsection{Laguerre tessellations in hyperbolic space}

Unlike ordinary Voronoi tessellations, Laguerre tessellations have no unique
canonical extension to hyperbolic space. Replacing squared Euclidean distance
in \eqref{eq:power} with squared hyperbolic distance leads to weighted
diagrams arising in optimal transport \cite{ShiWang2020}, but their cells need
not be geodesically convex. Another approach uses projective geometry and
Minkowski polarity. Related constructions occur in non-Euclidean Laguerre
geometry \cite{BobenkoLutzPottmannTechter2021} and in the theory of
decorated hyperbolic surfaces \cite{Lutz2023}. Minkowski convex-hull methods
also underlie the classical Epstein-Penner decomposition of cusped hyperbolic
manifolds \cite{EpsteinPenner1988}. A similar approach based on projective geometry and polarity has been used in
\cite{Sugihara} for constructing Laguerre tessellation on the sphere. In computational geometry, the
ordinary hyperbolic Voronoi diagram can be represented as Euclidean power diagram
restricted to the Klein model, with weighted variants discussed in
\cite{NielsenNock2010}. 

We adopt a Lorentzian construction in which the difference of two power
functions is linear in Minkowski coordinates. Consequently, cells are
intersections of hyperbolic halfspaces, and nonempty bisectors lie in
totally geodesic hyperplanes.
Moreover, the construction is suitable and convenient for studying totally geodesic sections and
allows generators on the ideal boundary, both of which are essential for our
probabilistic results.

We work in the hyperboloid model
\[
  \mathbb{H}^d:=\{u\in\RR^{d+1}:\loren{u}{u}=-1,\ u_{d+1}>0\},
  \qquad
  \loren{x}{y}:=\sum_{i=1}^d x_i y_i-x_{d+1}y_{d+1}.
\]
The ideal boundary $\partial_\infty\mathbb{H}^d$ consists of future lightlike
rays. We fix the representatives
\[
  \partial_\infty\mathbb{H}^d\simeq
  \{(\theta,1):\theta\in\SPH^{d-1}\},
\]
and write $\widehat{\mathbb{H}^d}:=\mathbb{H}^d\cup\partial_\infty\mathbb{H}^d$ for the
geometric compactification. With this normalization, for
$(v,h)\in\widehat{\mathbb{H}^d}\times(0,\infty)$ and $u\in\mathbb{H}^d$ define the
\emph{hyperbolic power function}
\begin{equation}\label{eq:HypPower}
  \pow_{\rm hyp}((v,h),u):=-\loren{hv}{u}.
\end{equation}
For an interior nucleus $v\in\mathbb{H}^d$, this becomes
$\pow_{\rm hyp}((v,h),u)=h\,\cosh( d_{\rm hyp}(v,u))$.
For an ideal nucleus $\widehat\theta=(\theta,1)\in\partial_\infty\mathbb{H}^d$, it satisfies
\[
  \log\pow_{\rm hyp}((\widehat\theta,h),u)
   =b_\theta(u)+\log h,
  \qquad b_\theta(u):=\log\bigl(-\loren{\widehat\theta}{u}\bigr),
\]
where $b_\theta$ is the Busemann function normalized to vanish at
$o=(0,\ldots,0,1)$. Thus the same power function incorporates both
hyperbolic distance and Busemann functions. For interior nuclei, taking logarithms yields an additively weighted diagram with cost $\log\cosh d_{\rm hyp}$, rather than the squared-distance cost commonly used in optimal transport.

For a countable set $A\subset\widehat{\mathbb{H}^d}\times(0,\infty)$, the Laguerre cell
of $(v,h)\in A$ is
\begin{equation}\label{eq:deterministic-laguerre-cell}
  C_{\rm hyp}((v,h),A)
  :=\bigl\{u\in\mathbb{H}^d:\ \pow_{\rm hyp}((v,h),u)
  \leq\pow_{\rm hyp}((v',h'),u)
  \text{ for all }(v',h')\in A\bigr\}.
\end{equation}
We call $v$ its \emph{nucleus} and denote by $\mathcal L_{\rm hyp}(A)$
the family of cells with nonempty interior, called the \emph{hyperbolic
Laguerre diagram}. We reserve \emph{tessellation} for a locally finite
family of closed, geodesically convex sets with nonempty, pairwise disjoint
interiors that cover $\mathbb{H}^d$ (see Definition~\ref{def:tessellation}).
If all nuclei are in $\mathbb{H}^d$ and have equal weights, the diagram reduces
to the ordinary hyperbolic Voronoi diagram.

\subsection{Main results}
Our first results give sufficient conditions under which a
countable collection of weighted points generates a tessellation, as well as
conditions ensuring boundedness of cells and normality of tessellation (see Theorem~\ref{thrm:deterministic normal}).
We then consider random generating sets. Let $\mu$ be a nonzero, locally
finite, nonnegative Borel measure on $(0,\infty)$ and let $\lambda>0$.
Write $\eta_{\mu,\lambda}$ for a Poisson point process on
$\mathbb{H}^d\times(0,\infty)$ with intensity measure
\[
  \Lambda_{\mu,\lambda}
  :=\lambda\bigl(\Vol_{\mathbb{H}^d}\otimes\mu\bigr),
\]
and put
$\mathcal L_{\rm hyp}(\mu,\lambda)
:=\mathcal L_{\rm hyp}(\eta_{\mu,\lambda})$.
We write $\mathcal L_{\rm hyp}(\mu)$ when the factor $\lambda$ has been
absorbed into $\mu$, and use a superscript to indicate the dimension.

\begin{theorem}\label{prop:existence}
Let $d\geq 2$, $\lambda>0$, and let $\mu$ be a nonzero, locally finite,
nonnegative Borel measure on $(0,\infty)$. Then
$\mathcal L_{\rm hyp}(\mu,\lambda)$ is almost surely a tessellation of
$\mathbb{H}^d$ if and only if
\begin{equation}\label{eq:ConditionMu}
  \int_{(0,T]}h^{-(d-1)}\,\mu(\mathrm dh)<\infty
  \qquad\text{for every }T>0.
\end{equation}
Whenever this condition holds, the tessellation almost surely has bounded
cells and is face-to-face and normal.
\end{theorem}

We call these tessellations \emph{hyperbolic Poisson-Laguerre
tessellations}. If $\mu(\mathrm dh)=f(h)\,\mathrm dh$, we also write
$\mathcal L_{\rm hyp}(f,\lambda)$. An important property of this class is
its stability under intersections with totally geodesic subspaces. The
following result is a hyperbolic analogue of
\cite[Theorem~4.1]{GusakovaInWoldeLuebke2025}.

\begin{theorem}\label{thrm:intersecdistrib}
Let $d\geq2$, $\lambda>0$, and $1\leq k\leq d-1$.
Suppose that $\mu$ is a nonzero, locally finite, nonnegative Borel measure
on $(0,\infty)$ such that
$\mathcal L^d_{\rm hyp}(\mu,\lambda)$ is almost surely a tessellation.
For any $k$-dimensional totally geodesic subspace $L_k\subset\mathbb{H}^d$,
up to isometry,
\[
  \mathcal L^d_{\rm hyp}(\mu,\lambda)\cap L_k
  \stackrel{d}{=}
  \mathcal L^k_{\rm hyp}(\widetilde\mu,\lambda),
\]
where
\[
  \widetilde\mu(\mathrm dh)
  =\left[
     \omega_{d-k}\,h^k
     \int_{(0,h)}(h^2-t^2)^{\frac{d-k}{2}-1}
                   t^{1-d}\,\mu(\mathrm dt)
   \right]\mathrm dh,
\]
and $\omega_m$ denotes the surface area of the unit sphere
$\SPH^{m-1}$.
\end{theorem}

One important special case of the above theorem is sectional Poisson-Voronoi tessellation. Let $\mathcal V_{\rm hyp}^d(\lambda)$ denote the
hyperbolic Poisson-Voronoi tessellation of intensity $\lambda$.
It corresponds to $\mu=\delta_1$, and the preceding theorem yields the
following analogue of \cite[Corollary~4.2]{GusakovaKabluchkoThaele2024}.

\begin{theorem}\label{thrm:Voronoiintersection}
Let $d\geq2$, $\lambda>0$, and $1\leq k\leq d-1$. For every
$k$-dimensional totally geodesic subspace $L_k\subset\mathbb{H}^d$, up to
isometry,
\[
  \mathcal V^d_{\rm hyp}(\lambda)\cap L_k
  \stackrel{d}{=}
  \mathcal L^k_{\rm hyp}
     \left(f_{k,\frac{d-k}{2}-1},\lambda\right),
\]
where, for $\beta>-1$,
\begin{equation}\label{eq:BetaDensity}
  f_{k,\beta}(h)
  :=c_{\beta}
     h^k(h^2-1)^\beta\,\mathbf1_{\{h>1\}},\qquad c_{\beta}:=\frac{2\pi^{\beta+1}}{\Gamma(\beta+1)}.
\end{equation}
\end{theorem}

We call $\mathcal L^k_{\rm hyp}(f_{k,\beta},\lambda)$ a
\emph{hyperbolic $\beta$-Voronoi tessellation}. For all $\beta>-1$,
it is almost surely a normal, face-to-face tessellation with bounded cells.

A second family of random tessellations  arises by considering intensity measures with the power-law densities
\[
  g_a(h):=h^a,\qquad a>d-2.
\]
By Theorem~\ref{prop:existence} and its one-dimensional analogue, Proposition~\ref{prop:existence-line},
$\mathcal L^d_{\rm hyp}(g_a,\lambda)$ is a tessellation with
bounded cells. Its law does not depend on $\lambda>0$, since multiplying all weights by a common factor leaves the diagram unchanged and the
height-intensity measure is homogeneous. We therefore write
$\mathcal L^d_{\rm hyp}(g_a)$ and call this family the
\emph{hyperbolic Gaussian-Voronoi tessellations}. The terminology is
motivated by the Euclidean Gaussian-Voronoi model, since under the logarithmic
change of variable $s=\log h$, the intensity density becomes
$e^{(a+1)s}\,\mathrm ds$.

This family also appears naturally in sections of ideal Poisson-Voronoi
tessellations. In the boundary normalization fixed above, let $\zeta$ be a
Poisson point process on $\partial_\infty\mathbb{H}^d\times(0,\infty)$ with
intensity measure
\[
  \Lambda_\infty
   :=\Vol_{\SPH^{d-1}}\otimes h^{d-2}\,\mathrm dh.
\]
The representation of the ideal Poisson-Voronoi tessellation in \cite{DAchilleEtAl2026} is equivalent, up to an
inessential deterministic rescaling of the weights, to the diagram
$\mathcal L_{\rm hyp}(\zeta)$. We use this Lorentzian formulation to
identify the distributions of its proper totally geodesic sections.

\begin{theorem}\label{thrm:idealsection}
Let $d\geq2$ and $1\leq k\leq d-1$. Let $\mathcal V_\infty^d$ be the
ideal Poisson-Voronoi tessellation of $\mathbb{H}^d$ and let
$L_k\subset\mathbb{H}^d$ be a $k$-dimensional totally geodesic subspace.
Then, up to isometry,
\[
  \mathcal V_\infty^d\cap L_k
  \stackrel{d}{=}\mathcal L^k_{\rm hyp}(g_{d-2}).
\]
\end{theorem}

In particular, although every cell of $\mathcal V_\infty^d$ is almost
surely unbounded, every proper totally geodesic section is a
Poisson-Laguerre tessellation with bounded cells. Indeed, the exponent
$d-2$ lies strictly above the existence threshold $k-2$ in every
dimension $k<d$. Thus hyperbolic Poisson-Laguerre tessellations provide a
class stable under totally geodesic sections and a common description of
ordinary, sectional and ideal Poisson-Voronoi tessellations.

\paragraph{\textbf{Organization of the paper.}} The rest of the paper is organized as follows. Section \ref{sec:preli} is a preliminary section, containing relevant background on hyperbolic spaces, projective geometry and tessellations. In Section \ref{sec:construction} we present a few alternative constructions of hyperbolic Laguerre diagram, provide a motivation for the choice of power function and give a proof of Theorem \ref{thrm:deterministic normal}. The proof of Theorem~\ref{prop:existence} is given in Section~\ref{conditions}.
The proofs of Theorems \ref{thrm:intersecdistrib}, \ref{thrm:Voronoiintersection} and \ref{thrm:idealsection} are given in Section~\ref{sec:sectional}. 

\section{Preliminaries}\label{sec:preli}

We begin by recalling some basic notions from hyperbolic and projective geometry needed for the construction of hyperbolic Laguerre diagrams. We conclude the section with a formal definition of a tessellation.

\subsection{Hyperbolic space and its geometric models}
Hyperbolic space $\mathbb{H}^d$, $d\ge 1$ admits two alternative descriptions, Riemannian and projective. From the Riemannian point of view, $\mathbb{H}^d$ for $d\ge 2$ is the unique, up to isometry, complete, simply connected $d$-dimensional Riemannian manifold of constant sectional curvature $-1$. A concrete realization is given by the hyperboloid model, in which the Riemannian metric is induced by the Lorentzian inner product. From the projective viewpoint, hyperbolic space is realized in the Beltrami-Klein model, in which geodesics are straight line segments and the metric is expressed in terms of cross-ratios. In what follows we will use both models: the hyperboloid model for defining Laguerre tessellations and Beltrami-Klein model for studying their sectional properties. For  detailed description of these models, we refer to \cite{Ratcc,Cannon,Bridson}.
 
We begin with the \emph{hyperboloid model}. For $x=(x_1,\ldots,x_{d+1})$ and $y=(y_1,\ldots,y_{d+1})$ in $\mathbb{R}^{d+1}$, the \emph{Lorentzian inner product} and the associated \emph{quadratic form} are defined by 
$$
\loren{x}{y} := x_1 y_1 + \cdots + x_dy_d - x_{d+1} y_{d+1}, \quad \text{and}  \quad \Vert x \Vert^2_L :=  \loren{x}{x}.
$$
Hyperbolic space is identified with the upper sheet of two-sheeted hyperboloid 
$$
\mathbb{H}^d := \{\, (x_1,\ldots,x_d,x_{d+1}) \in \mathbb{R}^{d+1} : \Vert x\Vert_L^2 = -1 \ \text{and}\ x_{d+1} > 0\, \}.
$$
The restriction of the Lorentzian inner product to the tangent spaces  of $\mathbb{H}^d$ defines a Riemannian metric of constant sectional curvature $-1$. The corresponding volume measure is denoted by $\Vol_{\mathbb{H}^d}$ and the hyperbolic distance between $x,y \in \mathbb{H}^d$ is given by
$$
d_{\rm hyp}(x,y) := \operatorname{arcosh} \bigl(-\loren{x}{y}\bigr).
$$
Moreover, we denote by $B(u,r)$ the geodesic ball centered at $u\in\mathbb{H}^d$ of radius $r>0$.

For $0\leq k\leq d$ by $k$-\emph{subspace} of $\mathbb{H}^d$ we understand the $k$-dimensional totally geodesic submanifold of $\mathbb{H}^d$. In the hyperboloid model, such subspaces are the non-empty intersections of $\mathbb{H}^d$ with $(k+1)$-dimensional linear subspaces of $\mathbb{R}^{d+1}$.

The \emph{Beltrami-Klein model} is obtained from hyperboloid model via \emph{gnomonic projection}
\begin{align*}
\Gnom\colon \HH^{d}&\to \BB^d\\
(x_1,\ldots,x_d,x_{d+1}) &\mapsto  \left( \frac{x_1}{x_{d+1}}, \ldots, \frac{x_d}{x_{d+1}} \right),
\end{align*}
 where $\mathbb{B}^d$ denotes the open Euclidean unit ball. This map is a bijection with inverse
$$
\Gnom^{-1}(v)
=
\frac{1}{\sqrt{1-\|v\|^2}}
(v_1,\ldots,v_d,1),\qquad v\in\BB^d.
$$
An important feature of the Beltrami-Klein model is that its $k$-dimensional hyperbolic subspaces are precisely the non-empty intersections of $\BB^d$ with $k$-dimensional affine subspaces of $\mathbb{R}^d$. This property reflects the underlying projective structure, which also motivates the name \emph{projective model}. To explain this connection, we recall some basic notions from projective geometry.

The $d$-dimensional real projective space $\mathbb{RP}^d$ is the set of all one-dimensional linear subspaces of $\RR^{d+1}$. Equivalently, its points are equivalence classes $[x]=[x_1 : \ldots : x_{d+1}]$ of nonzero vectors $x=(x_1,\ldots,x_{d+1})\in \RR^{d+1}\setminus\{0\}$, where $x,y\in \RR^{d+1}\setminus\{0\}$ are equivalent if and only if $y=\lambda x$ for some $\lambda \in \RR\setminus 0$. The notation $[x_1 : \ldots : x_{d+1}]$ represents a point in \emph{homogeneous coordinates}. 

For our purposes, it is convenient to work with \emph{affine coordinates}. Consider the subset $U_d=\{[x]\in\mathbb{RP}^d\colon x_{d+1}\neq0\}$. The map $\alpha\colon U_d\mapsto \RR^d$, given by
\begin{equation}\label{eq:AffineCoordinates}
\alpha([x_1: \ldots : x_d : x_{d+1}]) := \left( \frac{x_1}{x_{d+1}},\ldots,\frac{x_d}{x_{d+1}} \right)
\end{equation}
defines an \emph{affine chart}, identifying $U_d$ with $\RR^d$.
Its inverse is given by $(u_1,\ldots,u_d)\mapsto [u_1:\cdots:u_d:1]$. The complement $\mathbb{RP}^d\setminus U_d$, consisting of points with $x_{d+1}=0$, is called the \emph{hyperplane at infinity} (relative to this affine chart).

The \emph{absolute quadric} associated with the Lorentzian form is the projective hypersurface  $\mathcal{Q}_{\rm hyp}=\{\,[x] \in \RP^d: \loren{x}{x}=0\,\}$. The projective model of hyperbolic space consists of the projective points lying in the interior of the absolute quadric, namely
\[
\mathbb{H}^d=\{[x]\in\mathbb{RP}^d\colon \loren{x}{x}<0\}.
\]
The ideal boundary $\partial_{\infty}\mathbb{H}^d$ is identified with  $\mathcal{Q}_{\rm hyp}$. Since every point $[x]$, satisfying $\loren{x}{x}\leq 0$ belongs to $U_d$, the affine chart identifies $\mathbb{H}^d$ with $\mathbb{B}^d$ and its ideal boundary $\partial_\infty\mathbb{H}^d$ with $\mathbb{S}^{d-1}$. Under this identification, the passage from the hyperboloid model to the projective model coincides with the gnomonic projection introduced above.

\subsection{Ideal boundary, Busemann functions and horoballs}

The ideal boundary $\partial_\infty\mathbb{H}^d$ can also be described in terms of equivalence classes of geodesic rays or Busemann functions.
Fix $o=(0,\ldots,0,1)\in\mathbb{H}^d$ (in hyperboloid model) and identify $\partial_\infty\mathbb{H}^d$ with $\mathbb{S}^{d-1}$ as above. For  $\theta\in\mathbb{S}^{d-1}$, the unit-speed geodesic ray starting at $o$ and pointing towards $\theta$ is $\gamma_\theta(t)=(\sinh(t)\,\theta,\cosh(t))$, $t\geq 0$. The corresponding \emph{Busemann function}, normalized by
$b_\theta(o)=0$, is
\[
b_\theta(x)
:=\lim_{t\to\infty}
\bigl(d_{\mathrm{hyp}}(x,\gamma_\theta(t))-t\bigr),
\qquad x\in\mathbb{H}^d.
\]
These functions also describe the horofunction boundary of
hyperbolic space. Writing $\widehat{\theta}=(\theta,1)$, we have
\[
\langle\widehat{\theta},\widehat{\theta}\rangle_L=0,
\qquad
b_\theta(x)
=\log\bigl(-\langle \widehat{\theta}, x\rangle_L\bigr).
\]
Here $\widehat{\theta}$ is a future-directed lightlike vector.

For $s\in\mathbb{R}$, the (closed) \emph{horoball} centered at
$\theta$ is the sublevel set
\[
B_{\theta,s}
:=\{x\in\mathbb{H}^d:b_\theta(x)\leq s\}
=\{x\in\mathbb{H}^d:
\langle \widehat{\theta}, x\rangle_L\geq -e^s\}.
\]
Its boundary is the \emph{horosphere} $\partial B_{\theta,s}$. Equivalently, every horoball is of the form
\[
B_{u,a}
=\{x\in\mathbb{H}^d:\langle x,u\rangle_L\geq-a\},
\]
where $a>0$ and $u\in\mathbb{R}^{d+1}\setminus\{0\}$ is
future-directed and lightlike. Its ideal center is the projective
point $[u]\in\partial_\infty\mathbb{H}^d$.

Thus horospheres arise by intersecting the hyperboloid with affine
hyperplanes of lightlike Lorentzian normal, and horoballs with the
corresponding affine half-spaces. This linear description in the ambient
Lorentzian space is useful in the construction of hyperbolic Laguerre
tessellations.
 
\subsection{Tessellations}\label{Tessellations}
We begin by recalling the deterministic definition of a tessellation of hyperbolic space, following \cite[Definition~2.3.14]{Herold}.

\begin{definition}\label{def:tessellation}
A collection $\mathcal{L}=\{C_i:i\in I\}$ of subsets of $\mathbb{H}^d$ is called a \emph{tessellation} of $\mathbb{H}^d$ if the following conditions hold
\begin{enumerate}
    \item Every cell $C_i$ is closed, geodesically convex, and has nonempty interior.
    \item The cells cover hyperbolic space,
    $\bigcup_{i\in I}C_i=\mathbb{H}^d$.
    \item The collection $\mathcal{L}$ is locally finite, that is, every compact subset of $\mathbb{H}^d$ intersects only finitely many cells.
    \item Distinct cells have disjoint interiors,
    $\operatorname{int}(C_i)\cap\operatorname{int}(C_j)=\varnothing$ for every $i\neq j$.
\end{enumerate}
\end{definition}

In contrast to the Euclidean setting \cite[Chapter~10]{SchneiderWeil}, cells are not required to be bounded. Every cell of a tessellation is necessarily a locally polyhedral convex set. A subset $C\subseteq\mathbb H^d$ is called \emph{locally polyhedral} at $v\in\mathbb H^d$ if there exist an open neighborhood $U$ of $v$ and finitely many closed geodesic half-spaces $H_1,\ldots,H_m$ such that
\[
C\cap U=U\cap H_1\cap\ldots\cap H_m.
\]  
The set $C$ is called locally polyhedral if it is locally polyhedral at any $v\in\mathbb H^d$.
In particular, cells may have infinitely many faces. At the same time, every bounded cell is a convex hyperbolic polytope. These properties follow from the corresponding Euclidean result \cite[Lemma 10.1.1]{SchneiderWeil} by passing to the Beltrami-Klein model; see also \cite{DAchilleEtAl2026}. 
By a random tessellation of $\mathbb{H}^d$ we understand a collection of closed subsets of $\mathbb{H}^d$ which is almost surely a tessellation (see \cite[Lemma 2.3.15]{Herold} for a measurability statement for further details).

A tessellation $\mathcal{L}$ is called \emph{face-to-face} if the intersection of any two distinct cells is either empty or a common face of both cells.  Given a face-to-face tessellation $\mathcal{L}=\{C_i:i\in I\}$ and a $k$-dimensional totally geodesic subspace $K\subset\mathbb{H}^d$, we write, with a slight abuse of notation,
\begin{equation}\label{def:intersec}
\mathcal{L}\cap K
:=\bigl\{C_i\cap K:i\in I,\ 
\operatorname{int}_K(C_i\cap K)\neq\varnothing\bigr\},
\end{equation}
where $\operatorname{int}_K$ denotes interior relative to $K$. The condition on the relative interior excludes intersections of cells with $K$ that have lower dimension and therefore are not cells of the induced tessellation. We are only interested in the case where $K$ will be a totally geodesic subspace of the space.

For $0\leq s\leq d-1$, an \emph{$s$-face} of $\mathcal{L}$ is an $s$-dimensional face of one of its cells. A face-to-face tessellation $\mathcal{L}$ is called \emph{normal} if every $s$-face of $\mathcal{L}$ is contained in exactly $d-s+1$ cells, for each $s\in\{0,\ldots,d-1\}$.

\section{Hyperbolic Laguerre diagram}\label{sec:construction}
 In this section, we develop the notion of a hyperbolic Laguerre diagram introduced in the introduction. As in the Euclidean setting \cite{GusakovaInWoldeLuebke2025}, we interpret Laguerre diagrams in terms of a crystallization process and present an alternative geometric construction based on projective geometry, polarity, and intersection sets. In particular, we show that the choice of power function in \eqref{eq:HypPower} is essentially unique in yielding convex (locally) polyhedral cells. Finally, we establish conditions on the set of generating points under which the resulting Laguerre diagram forms a tessellation in the sense of Definition~\ref{def:tessellation}.
 
\subsection{Laguerre diagram via crystallization process} 
If all generating nuclei belong to $\mathbb{H}^d$, i.e. $A\subset \mathbb{H}^d\times (0,\infty)$, the Laguerre diagram $\mathcal L_{\rm hyp}(A)$ can be interpreted as a crystallization process, where each weighted point $(v,h)$ represents a crystal nucleus $v$ with a prescribed birth time $h$. All crystals follow the same isotropic growth law, and each point of space is assigned to the crystal that reaches it first. The arrival time of the crystal $(v,h)$ at point $u\in\mathbb{H}^d$ is given by the power function
\[
\pow_{\rm hyp}((v,h),u)=-\loren{hv}{u}=h\,\cosh(d_{\rm hyp}(v,u)).
\]
Thus, at time $t$ the crystal process $(v,h)$ has reached  reached the horoball
$B_{\theta,\log(t/h)}$. The boundary between two neighboring crystals (if not empty) is determined by equal arrival times. Note that we may include ideal generators $(\theta, h)$ in our construction, but $h$ would not represent the birth time of such crystal. An ideal crystal $(\theta, h)$ is formally born at time $-\infty$ and at time $t$ it has reached the set of points $B_{\theta, s+\log h}$. To avoid unnecessary complexity, in the rest of subsection we restrict to nuclei in $\mathbb{H}^d$.

Next we consider a more general model. Let $X$ denote either $\RR^d$ or $\mathbb{H}^d$, equipped with the corresponding distance $d_X$. Consider a crystallization process described above, where the arrival time of the crystal process $(v,r)$, $v\in X$, $r\in \mathbb{R}$ at point $u\in X$ is given by $r+D(d_X(u,v))$ for some continuous, strictly increasing, and unbounded function $D:[0,\infty)\to[0,\infty)$ satisfying $D(0)=0$. Thus, given two crystals $(v,r)\neq(v',r')$ the set
\[
H((v,r),(v',r')):=\{u\in X\colon D(d_X(u,v))-D(d_X(u,v'))=r'-r\}
\]
describes all points which are reached by both crystals simultaneously in time. For our purposes we look for the function $D$, such that $H((v,r),(v',r'))$ is either empty or a totally geodesic subspace of dimension $d-1$ for any two weighted points $(v,r)$ and $(v',r')$.

\begin{proposition}\label{prop:D}
Let $d\ge 2$ and let $D:[0,\infty)\to[0,\infty)$ be strictly increasing and continuous with $D(0)=0$. Assume that for any distinct $(v,r),(v',r')\in X\times \RR$ the set $H((v,r),(v',r'))$ is either empty or a totally geodesic hyperplane. Then
\[
D(s)=\begin{cases}
    cs^2,\qquad &X=\mathbb{R}^d,\\
    c\log(\cosh s),\qquad &X=\mathbb{H}^d,
\end{cases},\qquad c>0.
\]
\end{proposition}

\begin{proof}
Let $v,v'\in X$ be distinct points and denote by $t:=d_X(v,v')>0$. Choose $r,r'\in \mathbb{R}$ such that $r'-r=-D(t)$. Then the set $H((v,r),(v',r'))$ is nonempty, since it contains a point $v$, and, hence, by assumption is a totally geodesic hyperplane. Further note that $H((v,r),(v',r'))$ is invariant under all isometries fixing $v$ and $v'$. Since it contains $v$ but not $v'$, rotational symmetry implies that it must be the hyperplane through $v$ perpendicular to the geodesic joining $v$ and $v'$.

Now choose any point $u\in H((v,r),(v',r'))$ and let $s=d_X(u,v)$. Then 
\begin{equation}\label{eq:Dequation}
    D(d_X(u,v'))=D(t)+D(s).
\end{equation}

For $X=\mathbb{R}^d$ the Pythagorean theorem in combination with \eqref{eq:Dequation} leads to equation
\[
D(\sqrt{t^2+s^2})=D(t)+D(s),\qquad t,s\ge 0.
\]
With $F(x)=D(\sqrt{x})$ and substitution $x=t^2$, $y=s^2$ the above equation is equivalent to additive Cauchy equation $F(x+y)=F(x)+F(y)$, $x,y\ge 0$, whose continuous solutions are $F(x)=cx$.  

For $X=\mathbb{H}^d$ the Pythagorean theorem and \eqref{eq:Dequation} yields
\[
D({\rm arcosh}(\cosh s\cosh t))=D(s)+D(t),\qquad s,t\ge 0.
\]
By taking $F(x)=D({\rm arcosh} x)$ and substituting $x=\cosh s$, $y=\cosh t$ we obtain multiplicative Cauchy equation $F(xy)=F(x)+F(y)$, $x,y\ge 1$. Its continuous solutions are of the form $F(x)=c\log x$.

In both cases strict monotonicity implies $c>0$ which finishes the proof.
\end{proof}

Note, that by Proposition \ref{prop:D} in hyperbolic case the resulting arrival time at $u$ is $r+c\log\cosh (d_{\rm hyp}(v,u))$.
The crystal originating from $(v,r)$ reaches $u$ no later than the
crystal originating from $(v',r')$ precisely when
$$
 c\log (\cosh(d_{\rm hyp}(v,u)))+r
\leq
c \log (\cosh(d_{\rm hyp}(v',u)))+r'.
$$
Since exponentiation preserves order the above is
equivalent to
$$
e^r\cosh (d_{\rm hyp}(v,u))^c
\leq
e^{r'}\cosh (d_{\rm hyp}(v',u))^c.
$$
This could be the power function for hyperbolic space, but to simplify the calculations that follow, we introduce the positive weights
$h:=e^{r/c}$ and $h':=e^{r'/c}$ and use the power function $\pow_{\rm hyp}((v,h),u)=h\,\cosh( d_{\rm hyp}(v,u))$. The same procedure in Euclidean space leads to the choice $\pow_{\rm Euc}((v,r),u)=\| v-u\|^2+r$.

\begin{remark}
For distinct nuclei, write $\delta=d_{\rm hyp}(v,v')>0$.
The ratio
$\cosh (d_{\rm hyp}(v,u))/\cosh (d_{\rm hyp}(v',u))$ ranges over
$(e^{-\delta},e^\delta)$. Its level sets are totally geodesic
hyperplanes perpendicular to the geodesic through the nuclei.
Hence $H((v,r),(v',r'))$ is non-empty precisely when $|r-r'|<d_{\rm hyp}(v,v')$.
If this condition fails, one crystal reaches every point of
$\mathbb{H}^d$ strictly before the other. In contrast, two distinct Euclidean
nuclei $v\neq v'$ always have a non-empty simultaneous arrival hyperplane $H((v,r),(v',r'))$
regardless of their birth times $r,r'$.
\end{remark}

\begin{remark}
In dimension $d=2$, the positive time parametrization gives linear
area growth before collisions. Indeed, the area of a hyperbolic
ball of radius $r$ is $2\pi(\cosh r-1)$.
For a weighted point $(v,h)\in\HH^2\times(0,\infty)$ and a time $s\geq h$, the crystal has radius $\operatorname{arcosh}\left(\frac{s}{h}\right)$. Its area is therefore $2\pi\left(\frac{s}{h}-1\right)$.
Thus the area of each crystal grows linearly in the positive time parameter $s$, with growth rate $2\pi/h$.
This is analogous to the Euclidean planar construction, where the choice $D(s)=s^2$ also makes the area of a growing crystal linear in time.
\end{remark}

\subsection{Laguerre diagram as projection of the polar intersection set}\label{projective}

In this subsection we provide an alternative geometric construction of hyperbolic Laguerre diagram. Recall that $\widehat{\mathbb{H}^d}=\mathbb{H}^d\cup\partial_\infty\mathbb{H}^d$. Consider a countable set of weighted points $\widehat{\mathbb{H}^d} \times(0,\infty)$ and corresponding Laguerre diagram $\mathcal L_{\rm hyp}(A)$. The main idea is to embed the set of generating points $A$ into the projective space $\RP^{d+1}$ and to associate to every weighted point a polar hyperplane with respect to suitable quadric. The obtained set of hyperplanes determine an intersection set in projective space. Projecting the facets of this set back onto $\mathbb{H}^d$ we obtain a subdivision of $\mathbb{H}^d$, which is precisely the hyperbolic Laguerre diagram $\mathcal L_{\rm hyp}(A)$. 

 In this section we work with hyperboloid model of hyperbolic space and fix representatives for ideal boundary $\partial_\infty\mathbb{H}^d\simeq
  \{(\theta,1):\theta\in\SPH^{d-1}\}$. To avoid confusion, we specify the notation used throughout this subsection. We use $v,u,w$ for points of $\mathbb{H}^d$ in hyperboloid model representation;  $x,y$ for arbitrary vectors in $\RR^{d+1}$; and $\tilde{x}$ for vectors in $\RR^{d+2}$.

We will introduce an additional projective coordinate to encode the weight $h\in(0,\infty)$. Since $\widehat{\mathbb{H}^d}$ can be realized within $\RP^d$ (see Section \ref{sec:preli}), we embed the weighted space $\widehat{\mathbb{H}^d}\times(0,\infty)$ into $\RP^{d+1}$ using the map
\begin{align*}
    \phi: \widehat{\mathbb{H}^d}\times(0,\infty)&\to \RP^{d+1}\\
    (v,h)&\mapsto [v_1:\cdots:v_{d+1}:1/h].
\end{align*}
 
We now describe the image of this embedding in affine coordinates. Let $[\xx]$ be a point of $\RP^{d+1}$ with $[\xx]\in U_{d+1}$ and affine coordinates $\alpha[\xx]$ defined by \eqref{eq:AffineCoordinates}. If $[\xx]=[v_1:\ldots:v_{d+1}:1/h]$ is the projective point corresponding to $(v,h)\in\widehat{\mathbb{H}^d}\times(0,\infty)$, then $[\xx]\in U_{d+1}$ and $\alpha[\xx]=hv$. This motivates our choice of $1/h$ as the last homogeneous coordinate, since $\pow_{\rm hyp}((v,h),u) = -\loren{hv}{u}$. Moreover, if $v\in \mathbb{H}^d$, then $\|v\|_L^2=-1$ and $v_{d+1}>0$, leading to $ \left\|\alpha[\xx]\right\|_L^2 = -h^2<0$ and $hv_{d+1}>0$. Similarly, if $v\in \partial_{\infty}\mathbb{H}^d$, then $\|v\|_L^2=0$ and $v_{d+1}=1$, implying $\loren{\alpha[\xx]}{\alpha[\xx]}=0$ and $hv_{d+1}>0$.

We denote the future time-like cone as 
\[
\mathcal{C}_{\mathrm{time}} := \left\{ x\in\RR^{d+1}: \|x\|^2_L< 0 \text{ and } x_{d+1}>0 \right\},
\]
the future light-like cone as 
\[
\mathcal{C}_{\mathrm{light}}=\left\{ x\in \RR^{d+1}: \| x\|^2_L=0 \text{ and } x_{d+1}>0\right\},
\]
and set $\mathcal{C}:=\mathcal{C}_{\mathrm{light}}\cup \mathcal{C}_{\mathrm{time}}$. We will show that $\widehat{\alpha}(\mathbb{H}^d\times (0,\infty))=\mathcal{C}_{\mathrm{time}}$ and  $\widehat{\alpha}(\partial_{\infty}\mathbb{H}^d\times (0,\infty))=\mathcal{C}_{\mathrm{light}}$, where 
\[
    \widehat{\alpha}:=\alpha\circ\phi :
    (v,h)\mapsto hv.
\]
Indeed, the inclusions $\widehat{\alpha}(\mathbb{H}^d\times (0,\infty))\subset\mathcal{C}_{\mathrm{time}}$ and $\widehat{\alpha}(\partial_{\infty}\mathbb{H}^d\times (0,\infty))\subset\mathcal{C}_{\mathrm{light}}$ are obvious by the above argument. At the same time $\widehat{\alpha}$ is a bijection, which we show by constructing the inverse map $\widehat{\alpha}^{-1}$ from $\mathcal{C}$ to $\widehat{\mathbb{H}^d}\times(0,\infty)$. For every $x\in\mathcal{C}$, define $H(x):=\sqrt{-\|x\|^2_L}$. Further for $x\in\mathcal{C}_{\mathrm{time}}$ set  \begin{equation}\label{eq:projectionHd}
    \proj_{\mathbb{H}^d}(x)
:= x / \sqrt{-\|x\|^2_L}\in \mathbb{H}^d,
\end{equation} 
and for $x\in\mathcal{C}_{\mathrm{light}}$ define  $\proj_{\partial_{\infty}\mathbb{H}^d}(x):= x/x_{d+1}\in \partial_{\infty}\mathbb{H}^d$.
The map $\widehat{\alpha}^{-1}$ can then be described as 
\begin{equation}\label{def:phi}
\widehat{\alpha}^{-1} (x) =
\begin{cases}
\left(
\proj_{\mathbb{H}^d}(x),H(x)
\right), \quad & \text{if } x \in  \mathcal{C}_{\mathrm{time}} \\ 
(\proj_{\partial_{\infty}\mathbb{H}^d}(x),x_{d+1}), & \text{if } x \in  \mathcal{C}_{\mathrm{light}}.
 \end{cases}
\end{equation}

\begin{remark}
For each $s>0$, define the level set $\Sigma(s) := \left\{ x\in\RR^{d+1}: \|x\|^2_L=-s^2 \text{ and } x_{d+1}>0 \right\}$. Thus, it holds that $\widehat{\alpha}(\mathbb{H}^d\times\{s\})=\Sigma(s)$ and $\mathcal{C}_{\mathrm{time}} = \bigcup_{s>0} \Sigma(s)$.  
Likewise each $s>0$, defines a level set $\Sigma^{\infty}(s):=\{x\in\RR^{d+1}: \| x\|_L^2=0 \text{ and } x_{d+1}=s\}$, for which it holds $\widehat{\alpha}(\partial_{\infty}\mathbb{H}^d\times\{s\})=\Sigma^{\infty}(s)$ and $\mathcal{C}_{\mathrm{light}}=\cup_{s>0}\Sigma^{\infty}(s)$. Note that even though $\mathcal{C}_{\mathrm{light}}\cup\{0\}$ is the boundary of $\mathcal{C}_{\mathrm{time}}$ in $\RR^{d+1}$, this is not true for the corresponding  level sets. Each point on $\mathcal{C}_{\mathrm{light}}$ is a limiting object. 
For a sequence of points in $\mathcal{C}_{\mathrm{time}}$ to converge to a point in $\mathcal{C}_{\mathrm{light}}$, the weights must converge to $0$. Recall that, in the crystal interpretation, a crystal growing from a boundary point can reach a point in the interior only if its growth began at the beginning of time.
\end{remark}

Using map $\widehat{\alpha}^{-1}$ allows us to carry out the remainder of the construction in $\RR^{d+1}$
and avoid the terminology of $\RP^{d+1}$, although the construction of the polar intersection set originates from the projective correspondence described above. The
reason for approaching the construction through projective geometry is
that, under the embedding above, $\mathbb{H}^d$ corresponds to the component
with $x_{d+1}>0$ of the affine part of a quadric hypersurface in
$\RP^{d+1}$. In affine coordinates, this component coincides with
$\mathbb{H}^d$ itself. We elaborate.

Let $\tilde{x}=(x_1,\ldots,x_{d+2})\in\RR^{d+2}$, and consider the quadratic form
$$
q(\tilde{x}):= x_1^2+\cdots+x_d^2-x_{d+1}^2+x_{d+2}^2=\|x\|^2_L+x_{d+2}^2.
$$
The corresponding quadric hypersurface in $\RP^{d+1}$ is $\mathcal{Q}:= \left\{ [\tilde{x}]\in\RP^{d+1}: q(\tilde{x})=0 \right\}$. For a point $[\tilde{x}]\in\mathcal{Q}$ with $x_{d+2}\neq0$, the equation $q(\tilde{x})=0$ is equivalent to $\Vert \alpha[\tilde{x}]\Vert^2_L=-1$. Thus, the affine part of the quadric in affine chart $x_{d+2}\neq 0$ is $\mathcal{Q}_{\rm aff}=\{ x\in\RR^{d+1}: \|x\|^2_L=-1\}$ and its component defined by $x_{d+1}>0$ is precisely $\mathbb{H}^d$. The part at infinity of the projective quadric $\mathcal{Q}$ is naturally homeomorphic to the ideal boundary $\partial_\infty\mathbb H^d$, both being identified with the sphere $\mathbb S^{d-1}$.

The $(d+2)\times (d+2)$ matrix associated with $q$ is $B:=\operatorname{diag} (1,\ldots,1  ,-1,1)$.
Since $B$ is symmetric and nonsingular, we can define a polarity relation on $\RP^{d+1}$ with respect to $\mathcal{Q}$. The polar hyperplane of a point $[\tilde{x}]\in\RP^{d+1}$ is
$[\tilde{x}]^\circ
:=
\left\{
[\tilde{y}]\in\RP^{d+1}:
\yy B\xx^T = 0
\right\}$. 
In the affine coordinates, the polar hyperplane of $x\in\mathcal{C}$ is represented as 
\[
x^\circ = \left\{ y\in\RR^{d+1}: \loren{x}{y}=-1 \right\},
\]
which is an affine hyperplane. If $x\in\Sigma(1)=\mathbb{H}^d$, then its polar hyperplane $x^\circ$ is the hyperplane tangent to $\mathbb{H}^d$ at $x$. More generally, suppose that $x\in\Sigma(1/h)$ for some $h>0$. Then $h x\in\Sigma(1)=\mathbb{H}^d$, and $x^\circ$ is parallel to the hyperplane tangent to $\mathbb{H}^d$ at $h x$. Moreover, $h^2 x\in\Sigma(h)$ belongs to $x^\circ$. Thus, $x^\circ$ is the translate of the hyperplane tangent to $\mathbb{H}^d$ at $h x$ that passes through $h^2 x$. If $x=s(\theta,1)\in\Sigma^{\infty}(s)$, $s>0$ then its polar hyperplane $x^\circ$ has normal vector $(-\theta,1)$ and distance $1/(\sqrt{2}s)$ from the origin.

 We now define the hyperbolic Laguerre diagram by radially projecting
the relevant faces of a convex set in $\mathbb{R}^{d+1}$ onto $\mathbb{H}^d$.
Let $\widehat{A}\subset\ctime$ be a countable non-empty set. For each $x\in \widehat{A}$, the polar
hyperplane
$x^\circ$ divides $\RR^{d+1}$ into two closed halfspaces. We denote by
$(x^\circ)^+$ the halfspace that does not contain the origin, namely
$(x^\circ)^+
:=
\left\{
y\in\RR^{d+1}:
\loren{y}{x}\leq-1
\right\}.$ The \emph{polar intersection set} generated by $\widehat{A}$ is defined by
$$
P(\widehat{A})
:=
\bigcap_{x\in \widehat{A}}(x^\circ)^+.
$$
The set $P(\widehat{A})$ is closed and convex. Whenever it is non-empty, it is
unbounded. We note that if $\widehat{A}$ is finite $P(\widehat{A})$ is a polyhedral set, but in general it is not the case. In next section we will show that under suitable conditions on $\widehat{A}$ the set $P(\widehat{A})$ is a locally polyhedral in $\RR^{d+1}$ at every point of $\mathcal{C}_{\mathrm{time}}$.

For each $x\in \widehat{A}$, define
$
F(x,\widehat{A})
:=
P(\widehat{A})\cap x^\circ$.
We call the set $F(x,\widehat{A})$ the (exposed) face of $P(\widehat{A})$ induced by $x^\circ$ and note that it may be
empty. If it has dimension $d$, then it is a facet of $P(\widehat{A})$. The
Laguerre cell corresponding to $x$ is defined by
$$
C_{P}(x,\widehat{A})
:=
\proj_{\mathbb{H}^d}
\left(
F(x,\widehat{A})\cap\mathcal{C}_{\mathrm{time}}
\right),
$$
where $\proj_{\mathbb{H}^d}$ is defined in \eqref{eq:projectionHd}.
Let $\mathcal{F}_d(P(\widehat{A}))$ denote the collection of $d$-dimensional
facets of $P(\widehat{A})$. The hyperbolic Laguerre diagram generated by $\widehat{A}$ is
defined by
$$
\mathcal{L}_P(\widehat{A})
:=
\left\{
\proj_{\mathbb{H}^d}\left(F\cap\mathcal{C}_{\mathrm{time}}\right):
F\in\mathcal{F}_d(P(\widehat{A})),
\ F\cap\mathcal{C}_{\mathrm{time}}\neq\emptyset
\right\}.
$$
The next proposition shows that the construction through the polar
intersection set agrees with the definition of the hyperbolic Laguerre diagram
through the power function.

\begin{proposition}\label{prop:polyhedral}
Let $\widehat{A}\subset\ctime$ be a non-empty set, and let $\widehat{\alpha}^{-1}(\widehat{A}) \subset \widehat{\mathbb{H}^d}\times(0,\infty)$ be the corresponding set of weighted points. Then, for every $x\in \widehat{A}$,
$$
C_{P}(x,\widehat{A})
=
C_{\rm hyp}\bigl(\widehat{\alpha}^{-1} (x),\widehat{\alpha}^{-1} (\widehat{A})\bigr).
$$
\end{proposition}
\begin{proof}
Fix $x\in \widehat{A}$ and $v\in\mathbb{H}^d$, and consider the ray $\{rv:r>0\}\subset\mathcal C_{\mathrm{time}}$. For every $y\in \widehat{A}$ we have $\loren{v}{y}<0$. Hence this ray intersects the polar hyperplane $y^\circ$ at a unique point $r_y(v)v$, where
\begin{equation}\label{eq:radial-parameter}
r_y(v)=-\frac{1}{\loren{v}{y}}.
\end{equation}
For $r>0$, we have that $ rv\in(y^\circ)^+$ if and only if $r\loren{v}{y}\leq-1$. Since $\loren{v}{y}<0$, this is equivalent to $r\geq r_y(v)$.

By the definition of radial projection, $v\in C_P(x,\widehat{A})$ if and only if the ray intersects $F(x,\widehat{A})=P(\widehat{A})\cap x^\circ$. Its unique intersection with $x^\circ$ is $r_x(v)v$, so $v\in C_P(x,\widehat{A})$ is equivalent to $r_x(v)v\in P(\widehat{A})$. By definition of $P(\widehat{A})$ as the intersection of all half-spaces $(y^\circ)^+$, $r_x(v)v\in P(\widehat{A})$ if and only if $r_x(v)v\in(y^\circ)^+$ for every $y\in \widehat{A}$. By the preceding observation, this is equivalent to $r_y(v)\leq r_x(v)$ for every $y\in \widehat{A}$. Substituting equation~\eqref{eq:radial-parameter}, and using the fact that $\loren{v}{x}<0$ and $\loren{v}{y}<0$, we obtain  
\begin{equation}\label{eq:powerinequalities}
\loren{v}{x} \geq \loren{v}{y}
\end{equation}
for every  $y\in \widehat{A}$.
 
For any $y\in\mathcal C$, write $\widehat{\alpha}^{-1}(y)=(u_y,h_y)$. The decomposition $y=h_yu_y$ holds for both boundary points and inside the space points. Consequently,
$$
\pow_{\rm hyp}\bigl(\widehat{\alpha}^{-1}(y),v\bigr) =-h_y\loren{v}{u_y} =-\loren{v}{y}.
$$
Thus the family of inequalities (\ref{eq:powerinequalities}) are precisely
$$
\pow_{\rm hyp}\bigl(\widehat{\alpha}^{-1}(x),v\bigr)
\leq
\pow_{\rm hyp}\bigl(\widehat{\alpha}^{-1}(y),v\bigr)
$$
for every $y\in \widehat{A}$, which by \eqref{eq:deterministic-laguerre-cell} implies $v\in C_{\rm hyp}\bigl(\widehat{\alpha}^{-1}(x),\widehat{\alpha}^{-1}(\widehat{A})\bigr)$. Since $v\in\mathbb{H}^d$ was arbitrary, the claimed equality follows.
\end{proof}
 
 \begin{remark}
The polar construction also recovers the crystallization process described previously. Let $x\in\ctime$ and write the decomposition $\widehat{\alpha}^{-1} (x)=(v,h)$, so that $x=hv$. At time $t\geq h$, the growth front of the crystal originating from $(v,h)$ is $\proj_{\mathbb{H}^d}\left(\Sigma(1/t)\cap x^\circ\right)$.
Indeed, if $z\in\Sigma(1/t)\cap x^\circ$ and
$u=\proj_{\mathbb{H}^d}(z)$, then $z=u/t$, and the equation
$\loren{z}{x}=-1$ is equivalent to $\pow_{\rm hyp}((v,h),u)=t$.
Consequently, the crystal at time $t$ is
$$
\bigcup_{h\leq s\leq t}
\proj_{\mathbb{H}^d}\left(\Sigma(1/s)\cap x^\circ\right)
=
\left\{
u\in\mathbb{H}^d:
\pow_{\rm hyp}((v,h),u)\leq t
\right\}.
$$
Here the distinction between boundary points and points inside the space becomes clear.
If $v\in\mathbb{H}^d$, then at time $t=h$, the set consists only of the point $v$, and it grows
from that point as $t$ increases. If $v\in \partial_{\infty}\mathbb{H}^d$, at any time $t=h$ the crystal is a horoball and it has already expanded to an infinite size. 
\end{remark}

\subsection{Sufficient conditions for hyperbolic Laguerre tessellations}\label{Random Laguerre tessellations in hyperbolic spaces}
A Laguerre diagram need not be a tessellation in the sense of Definition~\ref{def:tessellation}. This issue is not specific to hyperbolic geometry, it also arises in Euclidean space. Schlottmann~\cite{Schlottmann} and Lautensack~\cite{Laut1} introduced conditions under which a weighted point configuration in $\RR^d$ generates a Laguerre tessellation. Here, we formulate analogous sufficient conditions for a weighted point configuration $A\subset \widehat{\mathbb{H}^d} \times(0,\infty)$ that ensure that the hyperbolic Laguerre diagram $\mathcal{L}_{\rm hyp}(A)$, is a tessellation of $\mathbb{H}^d$. Additionally, we impose general position condition under which the tessellation is normal and a condition under which the cells are bounded.

\begin{con}\label{conditions p1-p4}
Let $d\geq 1$ and $A\subset \widehat{\mathbb{H}^d}\times (0,\infty)$,
\begin{enumerate} 
\item[(P1)]
	For every $u\in \mathbb{H}^d$ and every $b > 0$
	there exist only finitely many $(v, h) \in A$ such that
			\begin{align*}
				\pow_{\rm hyp}((v,h), u) \leq b.
    	 	\end{align*}  
\item[(P2)] 
	For every $\widehat{\theta}  \in \partial_{\infty}\mathbb{H}^d$ and every $a >0 $ 
	there exists $(v,h) \in A$ such that $v\in \mathbb{H}^d$ and
    		$$ -\loren{hv}{\widehat{\theta}}\leq a. $$
\item[(P3)]
    The set of vectors $\{hv\colon (v,h)\in A\}$ is in affine general position, namely any $m\leq d+2$ distinct vectors are affinely independent.
\end{enumerate}
\end{con} 

\begin{theorem}\label{thrm:deterministic normal}
Let $d\geq 1$ and $A\subset \widehat{\mathbb{H}^d} \times (0,\infty)$ be a non-empty countable set of weighted points.
\begin{enumerate}
\item If $A$ satisfies condition \textup{(P1)}, then $P(\widehat{\alpha}(A))$ is a locally polyhedral set at any $x\in \mathcal{C}_{\mathrm{time}}$ and the hyperbolic Laguerre diagram $\mathcal{L}_{\rm hyp}(A)$ forms a face-to-face tessellation of $\mathbb{H}^d$.
\item If $A$ satisfies condition~\textup{(P2)}, then every nonempty
Laguerre cell generated by $(v,h)\in A$ is bounded when
$v\in\mathbb H^d$ and unbounded, with its nucleus $v$ as its
unique end, when $v\in\partial_\infty\mathbb H^d$.
\item If, in addition to \textup{(P1)}, the set $A$ also satisfies \textup{(P3)}, then the hyperbolic Laguerre diagram $\mathcal{L}_{\rm hyp}(A)$ is a normal tessellation of $\mathbb{H}^d$. \label{as:2}
\end{enumerate}
\end{theorem}

\begin{remark}
Under assumptions \textup{(P1)}-\textup{(P3)}, the requirement of non-empty interior in the definition of $\mathcal L_{\rm hyp}(A)$ can be replaced by nonemptiness, as the following lemma shows. 
\end{remark}

\begin{lemma}\label{lem:full-dimensional-cells}Let $d\geq 1$ and $A\subset \widehat{\mathbb{H}^d} \times (0,\infty)$ be a non-empty countable set of weighted points that satisfies conditions~\textup{(P1)}-\textup{(P3)}.
Then every non-empty Laguerre cell has non-empty interior.
\end{lemma}

\begin{remark}
Let us point out that conditions (1) and (4) of Definition \ref{def:tessellation} hold without \textup{(P1)}. Condition \textup{(P1)} prevents the following obstruction. If the weights of a sequence of generators decrease sufficiently rapidly, their powers at the origin may tend to zero. Since every generator has strictly positive power there, the infimum is not attained, and the origin belongs to no cell, which would violate condition (2) of Definition \ref{def:tessellation}.  
\end{remark}

If we consider $A\subset \mathbb{H}^d\times(0,\infty)$, then condition \textup{(P2)} ensures that the cells are bounded.
In Euclidean space, the corresponding condition is much simpler. Namely, if $A\subset \RR^d \times \RR$, one requires that
$$\conv_{\RR^d}(\{v:(v,h)\in A\} ) =\RR^d,$$
where $\conv_{\RR^d}$ denotes the convex hull. It would be natural to introduce an analogous hyperbolic condition, however, this is not sufficient to guarantee that all cells are bounded. The corresponding example is constructed in the lemma below and can easily be transferred to any dimension.

\begin{lemma}\label{counterexample}
Let $\HH^1$ be the hyperbolic line. There exists a set of weighted points 
$A=\{ (v_n , h_n)   : n \geq 0 \}$ such that the hyperbolic convex hull satisfies $$\conv_{\mathbb{H}^d}(\{v_n:n\in\mathbb{N}\})=\HH^1,$$ but
$C_{\rm hyp}((v_0,h_0),A)=\HH^1$, and thus the corresponding cell is unbounded.
\end{lemma}

Condition~\textup{(P2)} has a natural horospherical interpretation.
For $\widehat{\theta}\in\partial_{\infty}\mathbb{H}^d$ and $a>0$, the set $\{
u\in\mathbb{H}^d:
-\loren{u}{\widehat{\theta}}\leq a \}$
is a horoball $B_{\theta,\log a}$ centered at the ideal point $\theta$. Condition~\textup{(P2)} requires that, for every
$\widehat{\theta}\in\partial_{\infty}\mathbb{H}^d$ and every $a>0$, there exists $(v,h)\in A$ such
that the nucleus $v$ belongs to the horoball $B_{\theta,\log a-\log h}$. Thus, the condition asks to find points with relevant weights arbitrarily deeply toward every point of the ideal boundary. 

In Euclidean space, the condition that the convex hull of the nuclei
is the entire space is equivalent to requiring that every halfspace
contain a nucleus. Condition~\textup{(P2)} is a weighted hyperbolic
analogue of this property, with horoballs 
playing the role of Euclidean halfspaces.

\begin{proof}[Proof of Theorem~\ref{thrm:deterministic normal}]
\textit{(1)} First we show that $\mathcal{L}_{\rm hyp}(A)$ is a tessellation, for which we will check conditions (1) -- (4) of Definition \ref{def:tessellation}. Given $(v,h),(v',h')\in A$ define the set 
\[
S((v,h),(v',h')):=\{u\in\mathbb{H}^d\colon \loren{u}{hv-h'v'}\geq 0\},
\]
which is a closed hyperbolic half-space, all of $\mathbb{H}^d$ or the empty set. The Laguerre cell of $(v,h)$ can then be described as
\[
C_{\rm hyp}((v,h),A)=\bigcap_{(v',h')\in A}S((v,h),(v',h')).
\]
Consequently every non-empty Laguerre cell is intersection of closed hyperbolic half-spaces with $\mathbb{H}^d$ and thus closed and geodesically convex, meaning that condition (1) is satisfied. Moreover, full dimensional cells have disjoint interiors justifying condition (4).

It remains to check conditions (2) and (3). For every $p=(v,h)\in A$, we denote $f_p(u):=-\loren{u}{hv}$, for every $u\in \mathbb{H}^d$.
Let $K\subset\mathbb{H}^d$ be a non-empty compact
set, then $K \subset B(o,R)$ for some $R>0$.
For $v\in \widehat{\mathbb{H}^d}$ we have  $ v_1^2+\cdots+v_d^2 \leq v_{d+1}^2 $. Writing $u\in\mathbb{H}^d$ in hyperbolic polar coordinates $u=(\phi \sinh r , \cosh r )$, where $\phi \in \SPH^{d-1}$, we observe that 
\begin{align*}
f_p(u)&= - h \big( \sinh r(\phi_1v_1+\ldots\phi_dv_d)  - \cosh r \, v_{d+1}
\big) \geq h v_{d+1} (\cosh r - \sinh r) = e^{-r} h v_{d+1} = e^{-r} f_p (o) .
\end{align*}
Therefore, $f_p(o) \leq e^R f_p(u)$ for every $u\in K$.

 Fix $p_0\in A$. Since $K$ is compact, the maximum $M:=\max_{u\in K} f_{p_0}(u)$ is obtained and we can define the set 
$$
\mathrm{rel}(K,A):=\{ p\in A: f_p(o) \leq e^R M +1 \} .
$$
Observe that any weighted point $p\in A$ satisfying $f_p(u) \leq f_{p_0}(u)$ for some $u\in K$ has to belong to $ \mathrm{rel}(K,A)$. Additionally by assumption \textup{(P1)}, the set $\mathrm{rel}(K,A)$ must be finite. Therefore for any point $u\in K$, the minimum of the power function is obtained by a point $p$ in $\mathrm{rel}(K,A)$ and, hence, $u\in C_{\rm hyp}(p,A)$. 
Although some cells may have empty interior, only finitely many
cells intersect $K$. Those with nonempty interior cover a dense
subset of $K$ and their union is closed and therefore covers all of $K$.
This implies that the diagram covers the space. Moreover, the covering is locally finite, since the number of cells hitting $K$ does not exceed the cardinality of $\mathrm{rel}(K,A)$. 

It remains to show that $\mathcal{L}_{\rm hyp}(A)$ is face-to-face. For any two cells $C_1:=C_{\rm hyp}((v_1,h_1),A)$ and $C_2:=C_{\rm hyp}((v_2,h_2),A)$ with $C_1\cap C_2\neq \emptyset$ it holds that
\[
C_1\cap C_2\subset\{u\in\mathbb{H}^d\colon \loren{u}{h_1v_1-h_2v_2}=0\}=:H,
\]
where $H$ is hyperbolic hyperplane. If $u\in C_1\cap H$, then $\loren{u}{h_1v_1}=\loren{u}{h_2v_2}\ge \loren{u}{h'v'}$ for any $(v',h')\in A$ and, hence, $u\in C_2$. This implies $C_1\cap C_2=C_1\cap H=C_2\cap H$. Since $H$ is supporting hyperplane of both $C_1$ and $C_2$ the claim follows.

The fact that $P(\widehat{\alpha}(A))$ is a locally polyhedral set at any $x\in \mathcal{C}_{\mathrm{time}}$ follows from Proposition \ref{prop:polyhedral} and the fact that $\mathcal{L}_{\rm hyp}(A)$ is a tessellation. Indeed for any open $U\subset \mathbb{H}^d$ with compact closure, there are only finitely many Laguerre cells intersecting $U$. Let $A'$ be the set of corresponding generating points and let $\widehat{U}=\{tv\colon v\in U,t>0\}$. Since cells cover the space we have
\[
P(\widehat{\alpha}(A))\cap \widehat{U}=\widehat{U}\cap \bigcap_{x\in\widehat{\alpha}(A')}(x^{\circ})^+,
\]
and the claim follows.

\medskip
\noindent
\textit{(2)} Let $C_{\rm hyp}((v,h),A)$ be a non-empty cell and let $u_n \in C_{\rm hyp}((v,h),A)$ be a sequence of points that converges to an ideal point $\widehat{\theta} \in \partial_{\infty}\mathbb{H}^d$. By definition then $u_n/ (u_n)_{d+1}$ converges to $\widehat{\theta}$ with the topology of $\RR^{d+1}$. Since $u_n\in C_{\rm hyp}((v,h),A)$, for every $(v',h')\in A$, 
$$-\loren{u_n}{hv}\leq -\loren{u_n}{h'v'}.
$$
We divide both sides by the $(d+1)$-th coordinate and take the limit to obtain
$$
-\loren{\widehat{\theta}}{hv}\leq -\loren{\widehat{\theta}}{h'v'}.
$$
By condition \textup{(P2)}, the right-hand side can be made arbitrary small. Since $-\loren{\widehat{\theta}}{hv}$ is non-negative, it follows that $\loren{\widehat{\theta}}{hv}=0$. If $v\in\mathbb{H}^d$, this is a contradiction, thus the cell is bounded. If $v\in\partial_{\infty}\mathbb{H}^d$, then equality implies that $v=\widehat{\theta}$, and no other ideal limit is possible. 

The cells that are non-empty and generated by a boundary  point have to be unbounded. Fix a point $u\in C_{\rm hyp}((v,h),A)$, where $v\in\partial_{\infty}\mathbb{H}^d$. Then for any other $(v',h')\in A$ we have 
\begin{equation}\label{eq:ineqbound}
-\loren{u}{hv}\leq -\loren{u}{h'v'} \quad \text{ and } \quad -\loren{v}{hv}\leq -\loren{v}{h'v'}.
\end{equation}
Therefore the geodesic ray from $u$ toward $v\in\partial_{\infty}\mathbb{H}^d$ satisfies inequality~(\ref{eq:ineqbound}). Hence the whole geodesic ray is contained in $C_{\rm hyp}((v,h),A)$.

\medskip
\noindent
\textit{(3)} First we note that by (1) the set $P(\widehat{\alpha}(A))$ is locally polyhedral in $\mathcal{C}_{\rm time}$. Moreover, as follows by the proof of Proposition \ref{prop:polyhedral} every $k$-dimensional face of $\mathcal{L}_{\rm hyp}(A)$ corresponds to $k$-dimensional face of $P(\widehat{\alpha}(A))$ intersecting $\mathcal{C}_{\rm time}$. Let $F$ be a $k$-dimensional face of $P(\widehat{\alpha}(A))$, and choose $y\in\operatorname{relint}(F\cap \mathcal{C}_{\rm time})$. Assume that exactly $m$ hyperplanes, corresponding to generators $(v_{i_1},h_{i_1}),\ldots,(v_{i_m},h_{i_m})$, contain $y$. By affine general position condition \textup{(P3)}, the vectors $h_{i_1}v_{i_1},\ldots,h_{i_m}v_{i_m}$ are affinely independent, implying $m\leq d+1$. Since they lie in the affine hyperplane $\{x\in\mathbb{R}^{d+1}:\langle y,x\rangle_L=-1\}$, which does not contain the origin, they are also linearly independent. Hence, since $P(\widehat{\alpha}(A))$ is locally polyhedral set in $\mathcal{C}_{\rm time}$, we conclude that $k=\dim F=d+1-m$. Consequently, $F$ is contained in exactly $m=d-k+1$ facets of $P(\widehat{\alpha}(A))$ and by Proposition \ref{prop:polyhedral} every $k$-dimensional face of $\mathcal{L}_{\rm hyp}(A)$  is contained in exactly $d-k+1$ cells. Thus the tessellation is normal.
\end{proof}

\begin{remark}\label{re:boundary}
In the proof of boundedness, we associated with every $u\in\mathbb{H}^d$ a normalized vector $\widehat{u} :=u/u_{d+1}$. This map can be viewed as the Beltrami--Klein projection, with its image written in the affine form
$$
\widehat{\BB^d}
:=\{ x\in\RR^{d+1} : x_{d+1}=1,\|x\|_L^2<0\},$$
whose boundary is $\partial_{\infty}\mathbb{H}^d$. However, this interpretation describes only the spatial information carried by the normalized vector.

Since $\widehat{\BB^d} \subset\mathcal{C}_\mathrm{time}$, each vector represents a weighted point. In particular, $\widehat u$ corresponds to a point with weight $1/u_{d+1}$. Thus, as $u$ escapes to infinity and $\widehat u$ approaches $\partial_{\infty}\mathbb{H}^d$, this weight tends to zero. The normalization therefore forces the weights to vanish as the nuclei approach the boundary. Although the starting times tend to zero, the relative timing information survives in the limiting description, and this is what the boundary weights encode. In this sense, all boundary points started at time zero, yet their weights retain the relative time differences between them.
\end{remark}

\begin{proof}[Proof of Lemma~\ref{lem:full-dimensional-cells}]
Let $C_{\rm hyp}((v_0,h_0),A)$ be non-empty and set $m:=\dim C_{\rm hyp}((v_0,h_0),A)$, meaning the dimension of the minimal totally geodesic subspace containing $C_{\rm hyp}((v_0,h_0),A)$. Let $\widehat A:=\widehat{\alpha}(A)$ and $x:=h_0v_0$. By Proposition \ref{prop:polyhedral} the cell $C_{\rm hyp}((v_0,h_0),A)$ is the radial projection of $F(x,\widehat A)\cap \mathcal C_{\mathrm{time}}$. Since radial projection preserves dimension and $C_{\mathrm{time}}$ is open it holds $\dim F(x,\widehat A)=m$.

Suppose, for a contradiction, that $m<d$. Then $F(x,\widehat A)$ is not a facet of $P(\widehat A)$. 
Choose $z\in F(x,\widehat A)\cap\mathcal C_{\mathrm{time}}$, which is possible
because $C_{\rm hyp}((v_0,h_0),A)$ is nonempty. By condition~\textup{(P1)},
in a neighbourhood of $z$ the polar intersection set $P(\widehat A)$ agrees with a
$(d+1)$-dimensional polyhedron $Q$ defined by finitely many of
the same polar inequalities, including the inequality corresponding to $x_0$. The face $Q\cap x_0^\circ$ has
dimension $m$, since it agrees locally with $F(x,\widehat A)$. It is therefore
contained in at least $d-m+1$ facets of $Q$.

Choose $d-m+1$ such facets, corresponding to weighted points $(v_1,h_1),\ldots,(v_{d-m+1},h_{d-m+1})\in A$. Since $m<d$ none of these facets is defined by $(v_0,h_0)$. We may choose affinely independent $z_0,\ldots,z_{m}\in F(x,\widehat A)\cap\mathcal C_{\mathrm{time}}$ sufficiently close to $z$, such that $\loren{z_j}{h_iv_i}=-1$, $0\leq i\leq d-m+1$, $0\leq j\leq m$. Since $\loren{z_j}{h_0v_0}=-1$ for all $0\leq j\leq m$ vectors  $z_0,\ldots,z_{m}$ are linearly independent and, hence, $h_0v_0,\ldots, h_{d-m+1}v_{d-m+1}$ belong to an $(d-m)$-dimensional affine subspace, meaning that they are affinely dependent. Since $d-m+2\leq d-2$ this contradicts~\textup{(P3)} and, thus, $m=d$.
\end{proof}

\begin{proof}[Proof of Lemma~\ref{counterexample}]
For $n\in\mathbb{N}$, set $v_n=(n,\sqrt{1+n^2})$, $u_n=(-n,\sqrt{1+n^2})$ and $h_n=\sqrt{1+n^2}+   n  $, note $(v_0,h_0)=(u_0,h_0)$. Let $A=\{ (v_n,h_n) : n \in\mathbb{N}\}\cup \{(u_n,h_n):n\in\mathbb{N} \}\cup \{(0,1)\}$.
The nuclei tend to both ideal endpoints as $n$ grows to infinity, so
their hyperbolic convex hull is \(\HH^1\).

For any $w=(x,y)\in\HH^1$, we have $|x|<y$, and hence for $n>0$
\begin{align*}
\pow_{\rm hyp}((v_n,h_n),w)
&=h_n\left( y\sqrt{1+n^2}  -x n \right)\\
&> h_n\left(\sqrt{1+n^2}- n \right) y\\
&= y
=\pow_{\rm hyp}((v_0,h_0),w).
\end{align*}
The same inequality holds for $u_n$, replacing $x$ with $-x$.
Reindexing the generators gives the configuration in the statement.
\end{proof}

\section{Poisson Laguerre tessellations}\label{conditions}

As previously discussed, we are interested in Laguerre diagrams generated by Poisson point processes. More precisely, in this section we consider the Poisson-Laguerre diagram
$\mathcal{L}_{\rm hyp}(\eta_{\mu,\lambda})$ generated by a Poisson
point process $\eta_{\mu,\lambda}$ on $\mathbb{H}^d\times(0,\infty)$
with intensity measure
\[
\Lambda_{\mu,\lambda}
:=\lambda\bigl(\Vol_{\mathbb{H}^d}\otimes\mu\bigr),
\]
where $\mu$ is a nonzero, locally finite Borel measure on
$(0,\infty)$ and $\lambda>0$. This intensity is homogeneous in
the spatial coordinate, while $\mu$ determines the intensity
in the weight coordinate. We will prove Theorem~\ref{prop:existence}, but first establish some lemmata. 
 
  \begin{proposition}\label{le:assumtotesselinside}
Let $d\ge 1$, and let $\mu$ be locally finite, non-trivial, Borel measure on $(0,\infty)$ and $\lambda>0$. Suppose that, for every $t>0$,
$$
\int_{(0,t]}
\left(
\int_0^{\frac{t^2}{h^2}-1}
\frac{z^{\frac d2-1}}{\sqrt{z+1}}\,\dd z
\right)\mu(\dd h)
<\infty.$$
Then the Poisson point process $\eta_{\mu,\lambda}$ on $\mathbb{H}^d$ almost surely satisfies \textup{(P1)}--\textup{(P3)}. Consequently,
$\mathcal{L}_{\rm hyp}(\mu,\lambda)$ is almost surely a normal, face-to-face tessellation of $\mathbb{H}^d$ with bounded cells.
\end{proposition}

\begin{remark}
    By using the change of variables $z=x/(1-x)$ we get
\[
\int_0^{\frac{t^2}{h^2}-1}
\frac{z^{{\frac{d}{2}}-1}}{\sqrt{1+z}}\,\mathrm dz
=\int_0^{1-\frac{h^2}{t^2}}
x^{{\frac{d}{2}}-1}(1-x)^{-\frac{d+1}{2}}\,\mathrm dx=B_{1-\frac{h^2}{t^2}}\Big({\frac{d}{2}},\frac{1-d}{2}\Big),
\]
where $B_x(a,b)$ denotes incomplete Beta function. For $d\ge 1$, the second parameter $(1-d)/2\leq 0$, but $B_x({\frac{d}{2}},\frac{1-d}{2})$ remains finite for $0\le x<1$, since the singularity around $1$ is excluded.
\end{remark}

\begin{proof}
\textit{Condition~\textup{(P1)}}. For every $t>0$ and $u\in\mathbb{H}^d$, define 
\begin{equation}\label{eq:Ft}
F_t(u):= \left\{ (v,h)\in \HH^d\times(0,\infty): \pow_{\rm hyp}((v,h),u) \leq t \right\}.
\end{equation}
We will show that  
\[
\mathbb{P}(\# \left( F_t(u)\cap \eta_{\mu,\lambda} \right) <\infty \text{ for all }u\in \mathbb{H}^d\text{ and }t\geq 0)=1,
\]
where $\#A$ denotes the cardinality of a set $A$. We start by computing
\begin{equation}\label{eq:expecnumpoi}
\Lambda_{\mu,\lambda} (F_t(u)) =
\lambda\int_{(0,\infty)}\int_{\mathbb{H}^d}
\mathbf{1}\left\{
\pow_{\rm hyp}((v,h),u) \leq t
\right\}
 \,\Vol_{\mathbb{H}^d}(\dd v)\, \mu( \dd h ).
\end{equation}
Since $\pow_{\rm hyp}((v,h),u)
=
h\,\cosh\left(d_{\rm hyp}(v,u)\right)$,
the condition $\pow_{\rm hyp}((v,h),u) \leq t$ is equivalent to
$
h \leq t$
and
$d_{\rm hyp}(v,u)
\leq \operatorname{arcosh}\left(t/h \right)$. The points satisfying the preceding condition have spatial coordinate $ v\in B\left( u,\operatorname{arcosh}\left(\frac{t}{h}\right) \right) $ and weight $h\in(0,t]$. Substituting into \eqref{eq:expecnumpoi}, we obtain
\begin{align}
\Lambda_{\mu,\lambda} (F_t(u)) 
&=
\lambda\int_{(0,t]}
\Vol_{\mathbb{H}^d}\left(
B\left(
u,\operatorname{arcosh}\left(\frac{t}{h}\right)
\right)
\right)
\, \mu ( \dd h ) .
\label{eq:lastballequat}
\end{align}
The volume of a hyperbolic ball does not depend on its center, thus
\eqref{eq:lastballequat} becomes
$$\Lambda_{\mu,\lambda} (F_t(u)) =
\omega_d\lambda
\int_{(0,t]}
\int_0^{\operatorname{arcosh}(t/h)}
\sinh^{d-1}(s)\,\dd s\, \mu (\dd h ).
$$
Set $z:=\sinh^2(s)$. Then
$\dd z
=
2\sinh(s)\cosh(s)\,\dd s
=
2\sqrt{z}\sqrt{z+1}\,\dd s$,
and hence
$\sinh^{d-1}(s)\,\dd s
=
\frac{1}{2}
\frac{z^{\frac d2-1}}{\sqrt{z+1}}\,\dd z$.
Moreover, $\sinh^2\left(\operatorname{arcosh}\left(\frac{t}{h}\right)\right)
=
\frac{t^2}{h^2}-1$, because $\sinh^2(y)=\cosh^2(y)-1$.
Therefore,
\begin{equation}\label{eq:LambdaMuFt}
\Lambda_{\mu,\lambda} (F_t(u))
=
\frac{\omega_d\lambda}{2}
\int_{(0,t]}
\int_0^{\frac{t^2}{h^2}-1}
\frac{z^{\frac d2-1}}{\sqrt{z+1}}\,\dd z\, \mu (\dd h)
<\infty
\end{equation}
by the assumption on $\mu$. Since $\Lambda_{\mu,\lambda}$ is the intensity measure of the Poisson point process, $\mathbb{P}(\# \left( F_t(u)\cap \eta_{\mu,\lambda}\right)<\infty )=1 $.
 
Let $D$ be a countable dense subset of $\mathbb{H}^d$ and the positive integers $\mathbb{Z}_{>0}$. Then the preceding calculations imply $\mathbb{P}(\cap_{(u,t)\in D\times \mathbb{Z}_{>0}} \{ \# \left( F_t(u)\cap \eta_{\mu,\lambda}\right)<\infty\})=1$. For every $u\in \mathbb{H}^d$ and $t\in (0,\infty)$ there is $u'\in D$ such that $d_{\rm hyp}(u,u')<1$. Since $\cosh(r+1)\leq e\cosh(r)$ we obtain $\pow_{\rm hyp}((v,h),u')\leq e \pow_{\rm hyp}((v,h),u)$. Choose $n\in\mathbb{Z}_{>0}$ such that $n> ec$ and observe that $ F_t(u)\subset F_n(u')$. Since $\# \left( F_n(u')\cap \eta_{\mu,\lambda}\right)$ is finite, so is $\#\left( F_t(u)\cap \eta_{\mu,\lambda} \right)$. Hence condition~\textup{(P1)} holds almost surely.

\textit{Condition~\textup{(P2)}}. Since $\mu$ is non-trivial, there exist $h_0<h_1$ such that $0<\mu([h_0,h_1])<\infty$.

In dimension $d\ge 2$ for each integer $n>2$, choose a maximal set $N_n$ on the hyperbolic sphere with center $o$ and radius $n$, whose distinct points have pairwise distance at least $1$. By maximality, every point of the sphere is within distance $1$ of a point of $N_n$. The open balls $B(x,1/2)$, $x\in N_n$, are disjoint and contained in $B(o,n+1/2)$, so $\# N_n  \leq \Vol_{\mathbb{H}^d}(B(o,n+1/2))/\Vol_{\mathbb{H}^d}(B(o,1/2))\leq Ce^{(d-1)n}$. 

Let $E_n$ be the event that, for some $x\in N_n$, the set $B(x,n/2-1)\times[h_0,h_1]$ contains no point of $\eta_{\mu,\lambda}$. By the union bound, $\mathbb P(E_n)\leq Ce^{(d-1)n}\exp\{-C'\Vol_{\mathbb{H}^d}(B(o,n/2-1))\}$ for some $C'>0$. This bound is summable in $n$, because the volume of balls grows exponentially. Thus, by Borel--Cantelli, almost surely, for all sufficiently large $n$, every such set contains a point of $\eta_{\mu,\lambda}$.

On this event, for any horoball take the geodesic ray from $o$ toward its ideal center. For all sufficiently large $n$, we have the ball of radius $n/2$ around the geodesic at height $n$ eventually is inside the horoball. By choosing a closest $z\in N_n$ we obtain that $B(z,n/2-1)$ is contained in the horoball and thus a point as well. Therefore, condition~\textup{(P2)} holds for all $b>0$ and all $\widehat{\theta}\in\partial_{\infty}\mathbb{H}^d$, almost surely.

In case $d=1$ we use another argument. Identify $\HH^1$ with $\RR$ by writing every point $v\in\HH^1$ as $v=( \sinh s, \cosh s)$ with $s\in\RR$. Almost surely, we can find nuclei with weights in $[h_0,h_1]$ with arbitrary large positive and negative $s$. Their powers at the ideal endpoints
$(1,1)$ and $(1,-1)$ are $he^{-s}$ and $he^s$, respectively. Since their weights are bounded, the infimum at each endpoint is zero, which gives \textup{(P2)}.

\textit{Condition \textup{(P3)}}. To verify \textup{(P3)}, first observe that under the affine map $\widehat{\alpha}$, the intensity measure assigns zero mass to any affine hyperplane of $\RR^{d+1}$. Indeed, a point $\widehat{\alpha}(v,h)=hv$ belongs to an affine hyperplane $\{x\in\mathbb{R}^{d+1}\colon\loren{z}{x}=-1\}$ if and only if $h\loren{z}{v}=-1$. For a fixed $h$ this condition means that $v$ belongs to a hyperbolic sphere, which has hyperbolic volume zero. The multivariate Mecke formula \cite[Chapter 4]{LastPenrose} therefore implies that condition \textup{(P3)} holds almost surely. 
\end{proof}

Since the measure $\mu$ is locally finite it can have a non-integrable singularity only around $0$. Thus, if the support of $\mu$ does not include $0$ the condition of Proposition \ref{le:assumtotesselinside} holds and we get the following corollary.

\begin{corollary}\label{cor:Conditions}
Let $d \geq 1$ and $\mu$ be a locally finite, non-trivial Borel measure on $(0,\infty)$ and $\lambda>0$. Suppose that $\operatorname{supp}( \mu)\subset [a,\infty)$, for some $a>0$. Then the Poisson point process $\eta_{\mu,\lambda}$ on $\mathbb{H}^d$ satisfies \textup{(P1)}-\textup{(P3)}, and therefore $\mathcal{L}_{\rm hyp}(\mu,\lambda)$ is almost surely a normal random tessellation of $\mathbb{H}^d$.
\end{corollary}

\begin{proof}
    Since for any $z\ge 0$ it holds that $\frac{z^{\frac d2-1}}{\sqrt{z+1}}\leq z^{\frac d2-1}$ we get for any $t\ge a$ that
    \[
    \int_{(0,t]}
\left(
\int_0^{\frac{t^2}{h^2}-1}
\frac{z^{\frac d2-1}}{\sqrt{z+1}}\,\dd z
\right)\mu(\dd h)\leq {\frac{2}{d}}\Big(\frac{t^2}{a^2}-1\Big)^{\frac{d}{2}}\mu([a,t])<\infty.
    \]
    For $t<a$ the above integral equals $0$ and the claim follows by Proposition \ref{le:assumtotesselinside}.
\end{proof}

Next we show that condition from Proposition \ref{le:assumtotesselinside} is equivalent to a simpler condition on the intensity measure $\mu$.

\begin{lemma}\label{le:equivalenceconditions}
 Let $d \geq 2$ and $\mu$ be a locally finite Borel measure. Then, for every $t>0$,
 $$
\int_{(0,t]} \int_{0}^{\frac{t^{2}}{h^{2}} - 1}\frac{z^{\frac{d}{2}-1}}{\sqrt{z+1}}\, \dd z\, \mu (\dd h) <\infty
\Leftrightarrow
\int_{(0,t]} h^{-(d-1)}\,\mu (\dd h)<\infty .
$$
\end{lemma}
\begin{proof}
Fix $t>0$. For $r\ge 1$, we set
$$J_d(r)
:=\int_0^{r^2-1}\frac{z^{\frac d2-1}}{\sqrt{z+1}}\,\dd z .$$
Then the double integral on statement of the lemma can be written as
$\int_{(0,t]} J_d(t/h)\,\mu( \dd h )$.

We first prove an upper bound. For every $z>0$ we have that $\sqrt{z+1}\ge \sqrt z$ and therefore
$$
\frac{z^{\frac d2-1}}{\sqrt{z+1}} \leq z^{\frac{d-3}{2}}.$$
Since we assumed that $d\ge2$ the exponent $(d-3)/2$ is strictly larger than $-1$ and from the integration we obtain the upper bound
$$
J_d(r)
\le
\int_0^{r^2-1}z^{\frac{d-3}{2}}\,\dd z
=
\frac{2}{d-1}(r^2-1)^{\frac{d-1}{2}}
\le
C_d r^{d-1},
$$
for some constant $C_d>0$ and all $r\ge1$.
The upper bound gives
$$
\int_{(0,t]}J_d(t/h) \,\mu(\dd h)
\leq C_d t^{d-1}
\int_{(0,t]}  h^{-(d-1)} \,\mu(\dd h),
$$
which proves one implication.

We will show that there is a constant $c_d>0$ such that $J_d(r)\geq c_d r^{d-1}$, for all $r\geq 2$. First consider $z\in[r^2/2,r^2-1]$, there  
$$\frac{z^{\frac d2-1}}{\sqrt{z+1}} \geq \frac{(r^2/2)^{\frac d2-1}}{r}
=2^{1-\frac d2}r^{d-3}.$$
So, integrating over $[r^2/2,r^2-1]$, we obtain
$$
\int_{\frac{r^2}{2}}^{r^2-1} \frac{z^{\frac d2-1}}{\sqrt{z+1}}\, \dd z \geq \left( \frac{r^2}{2}-1\right) 2^{1-\frac d2}r^{d-3} \geq 2^{-(1+\frac{d}{2})} r^{d-1}.
$$
Thus $J_d(r)\geq c_d r^{d-1}$, for all $r\geq 2$. 
Therefore, for $t/h \geq 2$, $h \in (0,t/2]$. 
This yields
$$
\int_{(0,t/2]}h^{-(d-1)} \, \mu(\dd h) \leq
\frac{1}{c_d t^{d-1}}
\int_{(0,t/2]}J_d(t/h) \, \mu(\dd h)
<\infty.
$$
Since $\mu$ is finite on $[t/2,t]$, we conclude,  
$\int_{(0,t]}h^{-(d-1)} \, \mu(\dd h) <\infty$,  and the desired equivalence is established. 
\end{proof}
\begin{proof}[Proof of Theorem~\ref{prop:existence}]
If the integrability condition holds, then
Proposition~\ref{le:assumtotesselinside} and Lemma~\ref{le:equivalenceconditions} show that
$\mathcal{L}_{\rm hyp}(\mu,\lambda)$ is almost surely a normal tessellation of $\mathbb{H}^d$.

Suppose that the integrability condition fails.
Then there exists
$t_0>0$ such that $\int_{(0,t_0]} h^{-(d-1)}\, \mu (\dd h ) =\infty$. Since $\mu$ is locally finite in $(0,\infty)$, the divergence can
only occur at the lower endpoint $0$. Thus, for every $t>0$ by Lemma \ref{le:equivalenceconditions} we get
$$
\int_{(0,t]} h^{-(d-1)}\,\mu (\dd  h )=\infty \Leftrightarrow  
 \int_{(0,t]} \int_{0}^{\frac{c^{2}}{h^{2}} - 1}\frac{z^{\frac{d}{2}-1}}{\sqrt{z+1}}\, \dd z\, \mu (\dd h) =\infty
$$

Fix now a point $u\in\mathbb{H}^d$ and for $t>0$ recall a definition \eqref{eq:Ft} of $F_t(u)$.
By \eqref{eq:LambdaMuFt} we have
$$\Lambda_{\mu,\lambda}(F_t(u)) =
\frac{\omega_d\lambda}{2}\int_{(0,t]} \int_0^{\frac{c^2}{h^2}-1}
\frac{z^{\frac d2-1}}{\sqrt{z+1}} \,\dd z\, \mu (\dd h )=\infty .
$$
Since $F_t(u)$ has infinite intensity measure, $\eta_{\mu,\lambda}$ has infinitely many points in $F_t(u)$ almost surely.
Applying this to $t_n=1/n$, we get that
$$\inf_{(v,h)\in\eta_{\mu,\lambda}}\pow_{\rm hyp}((v,h),u)\leq 0 .$$
But for any element $(v,h)\in\eta_{\mu,\lambda}$ the power function is $\pow_{\rm hyp}((v,h),u)
=
h\,\cosh\left(d_{\rm hyp}(v,u)\right)
\geq h>0$, therefore $u$ does not belong to any Laguerre cell.  
Thus the Laguerre diagram does not cover $\mathbb{H}^d$, and consequently it is not a tessellation.
\end{proof}
We now state the one dimensional analogue of Theorem~\ref{prop:existence}.
\begin{proposition}\label{prop:existence-line}
Let $\mu$ be a nonzero, locally finite Borel
measure on $(0,\infty)$ and $\lambda>0$. Then
$\mathcal L_{\rm hyp}^{1}(\mu,\lambda)$ is almost surely a
tessellation of $\HH^1$ if and only if
$$
\int_{(0,T]}\log(T/h)\,\mu(\dd h)<\infty
\qquad\text{for every }T>0.
$$
Whenever this condition holds, the tessellation almost surely is face-to-face, normal and
has bounded cells.
\end{proposition}
\begin{proof}
For one direction we note that for $T>0$
\[
\int_{(0,T]}
\left(
\int_0^{\frac{T^2}{h^2}-1}
\frac{z^{-1/2}}{\sqrt{z+1}}\,\dd z
\right)\mu(\dd h)=2\int_{(0,T]}\arcosh(T/h)\mu(\dd h).
\]
By local finiteness of $\mu$ and the bounds
$\log x\leq\operatorname{arcosh}(x)\leq\log x+\log 2$ for $x\geq1$,
this integral is finite if and only if the logarithmic integral in the
statement is finite. For the same reason, divergence for a $T>0$ implies divergence for every $T>0$. Hence, by Proposition \ref{le:assumtotesselinside} $\mathcal L_{\rm hyp}^{1}(\mu,\lambda)$ is almost surely a face-to-face, normal
tessellation of $\HH^1$ with only bounded cells.

The proof of other direction is analogous to the proof in Theorem \ref{prop:existence}.
\end{proof}

As a direct application of Corollary \ref{cor:Conditions} we can ensure that hyperbolic $\beta$-Voronoi tessellation $\mathcal L^k_{\rm hyp}(f_{k,\beta},\lambda)$ with $f_{k,\beta}$ defined by \eqref{eq:BetaDensity} is almost surely a normal face-to-face tessellation with bounded cells for $\beta>-1$. Also observe that according to Theorem~\ref{prop:existence} and Proposition~\ref{prop:existence-line} the hyperbolic Gaussian-Voronoi tessellation $\mathcal{L}_{\rm hyp} (g_a,\lambda)$ with $g_a(h)=h^a$ is a tessellation of $\HH^d$, if and only if $a > d-2$.

A special property of the hyperbolic Gaussian-Voronoi tessellation is that its distribution is invariant with respect to the intensity parameter. Given two different intensity measures with the same polynomial exponent, the corresponding inhomogeneous Poisson point processes define the same random tessellation object. The same invariance holds for the ideal
Poisson-Voronoi tessellation, as observed in \
\cite[Definition~2.1]{DAchilleEtAl2026}.
\begin{theorem}\label{thm:Gaussian}
For any $\lambda_1,\lambda_2>0$ and $a>d-2$ the Poisson-Laguerre tessellations $\mathcal{L}_{\rm hyp}(g_a,\lambda_1)$ and $\mathcal{L}_{\rm hyp}(g_a,\lambda_2)$ are equal in distribution.
\end{theorem}

\begin{proof}
Let $c > 0$ be a positive constant, and define the scaling mapping $T_c : \mathbb{H}^d \times (0,\infty) \to \mathbb{H}^d \times (0,\infty)$ by $T_c(v,h) = (v, ch)$. Observe that the Laguerre tessellation is invariant under $T_c$, namely $\mathcal{L}_{\rm hyp}(A) = \mathcal{L}_{\rm hyp}(T_c(A))$.

By the mapping theorem \cite[Theorem~24.16]{Klenke}, the point process $T_c(\eta_{g_a,\lambda_1})$ is a Poisson point process with intensity measure $\tilde{\Lambda}$ given by 
$$ 
\tilde{\Lambda}(A \times B)  
= \Lambda_{g_a,\lambda_1}(T_c^{-1}(A \times B)) 
= \lambda_1 \int_{\mathbb{H}^d} \int_{(0,\infty)}{\bf 1}\{ ch\in B, v\in A\} 
h^a \, \dd h \, \mathrm{Vol}_{\HH^d}(\dd v),
$$
for any measurable $A\subset \mathbb{H}^d$ and $B\subset (0,\infty)$. Applying the change of variables $w = ch$ we get
$$
\tilde{\Lambda}(A \times B) 
= \lambda_1 \int_A \int_B \left(\frac{w}{c}\right)^a c^{-1} 
\, \dd w \, \mathrm{Vol}_{\HH^d}(\dd v)  
= \lambda_1 c^{-(a+1)} \int_A \int_B h^a \, \dd h \, \mathrm{Vol}_{\HH^d}(\dd v).
$$
By choosing the scaling factor to be $c = (\lambda_1/\lambda_2)^{\frac{1}{a+1}}$ we have $\tilde{\Lambda}=\Lambda_{g_a,\lambda_2}$ and the claim follows.
\end{proof}

\section{Sectional properties of Laguerre tessellation}\label{sec:sectional}
In this section, we study intersections of Laguerre tessellations in $\mathbb{H}^d$ with a fixed totally geodesic subspace. Let $L_k$ be a $k$-dimensional subspace of $\mathbb{H}^d$. We show that the section $\mathcal{L}_{\rm hyp}^{d}(\mu)\cap L_k$, understood according to \eqref{def:intersec}, can be represented as a Poisson-Laguerre tessellation in $L_k\cong\HH^k$, and we determine the intensity measure of a Poisson point process that generates this random tessellation. Moreover, we determine the distribution of sectional Poisson-Voronoi and ideal Poisson-Voronoi tessellations.

Throughout this section we use a superscript to indicate the dimension of relevant objects (cones, tessellations, intersection sets) and subscript to indicate the dimension attributed to maps (i.e. $\widehat{\alpha}$, $\widehat{\alpha}^{-1}$). 

\subsection{Sections of Poisson-Laguerre tessellations}
The aim of this subsection is to prove Theorem~\ref{thrm:intersecdistrib}. We will split the proof in steps using a series of lemmata. 
For the preliminary lemmata, we restrict attention to
$$
\tilde{L}_k=\{(x_1,\ldots,x_{d+1}) : x_{k+1}=\ldots=x_d=0\}
$$
and $L_k = \mathbb{H}^d \cap \tilde{L}_k$. Since $L_k\cong\HH^k$ in what follows we will often identify $L_k$ with $\HH^k$. 
 Therefore, the cones $\mathcal C^k_{\mathrm{time}}$ and
$\mathcal C^k_{\mathrm{light}}$ are obtained by intersecting
$\mathcal C^d_{\mathrm{time}}$ and $\mathcal C^d_{\mathrm{light}}$,
respectively, with the linear span of $L_k$ in $\RR^{d+1}$.
Similarly, in affine coordinates, the lower-dimensional polarity
is induced by the intersection of the affine part of
$\mathcal Q^d$ with the same linear subspace. 

The strategy of the proof is based on the correspondence between the Laguerre tessellation $\mathcal{L}^d_{\rm hyp}(A)$ and locally polyhedral set $P^d(\widehat{\alpha}_d(A))$. Instead of studying the sectional tessellation $\mathcal{L}^d_{\rm hyp}(A)\cap L_k$ directly we will consider sectional locally polyhedral set $P^d(\widehat{\alpha}_d(A))\cap \tilde L_k$, whose facets determine the cells of sectional tessellation.

The first lemma shows that $P^d(\widehat{\alpha}_d(A))\cap \tilde L_k$ is the polar intersection set generated by the projections of points from $\widehat{\alpha}_d(A)$ onto $\tilde L_k$.

\begin{lemma}\label{prop:project}
Let $\widehat{A}\subset \ctime^d$, and let $\widehat{B}:=\proj_{\tilde{L}_k}(\widehat{A})\subset \ctime^k$, where
$$\proj_{\tilde{L}_k}(x_1,\ldots,x_{d+1}):=(x_1,\ldots,x_k,0,\ldots,0,x_{d+1})$$
is the linear projection onto $\tilde L_k$. Then
$$
P^d(\widehat{A})\cap\tilde{L}_k=P^k(\widehat{B}),
$$
where $P^k(\widehat{B})$ is the polar polyhedron generated by $\widehat{B}$ in
$\tilde{L}_k\cong \mathbb{R}^{k+1}$.
\end{lemma}

The next lemma expresses the projection operation directly in the
location-weight coordinates of $\widehat{\mathbb{H}^k}\times(0,\infty)$. For $v=(v_1,\ldots,v_{d+1})$, write
$\barbelow{v}:=(v_{k+1},\ldots,v_d)$, with $\|\barbelow{v}\|$
denoting its Euclidean norm.
\begin{lemma}\label{lem:coordexpre}
Let $\tau_k:=\widehat{\alpha}^{-1}_k\circ\proj_{\tilde L_k}\circ\,\widehat{\alpha}_d$ from $\widehat{\mathbb{H}^d}\times (0,\infty)$ to $\widehat{\HH^k}\times (0,\infty)$. Then we have
$$
\tau_k(v,h)
=
\begin{cases}
\left(
\frac{(v_1,\ldots,v_k,v_{d+1})}
     {\sqrt{1+\|\barbelow{v}\|^2}},
h\sqrt{1+\|\barbelow{v}\|^2}
\right), &v\in\mathbb{H}^d,\\
\left(
\frac{(v_1,\ldots,v_k,v_{d+1})}{\|\barbelow{v}\|},
h\|\barbelow{v}\|
\right),
& v\in\partial_{\infty}\mathbb{H}^d,\, \|\barbelow{v}\|>0,\\[6pt]
\bigl((v_1,\ldots,v_k,v_{d+1}),h\bigr),
& v\in\partial_{\infty}\mathbb{H}^d,\,\|\barbelow{v}\|=0.
\end{cases}
$$
\end{lemma}

The map $\tau_k$ first sends a weighted point in $\widehat{\mathbb{H}^d}\times(0,\infty)$ through $\widehat{\alpha}$ to its affine coordinates in $\ctime^{d}$. The projection $\proj_{\tilde L_k}$ then maps this point into the cone $\ctime^{k}$ contained in $\tilde L_k$. Finally, $\widehat{\alpha}^{-1}_k$ maps the projected point from $\ctime^{k}$ to its location-weight coordinates in $\widehat{\HH^k}\times(0,\infty)$. 

Having established the projection map, we now determine the intensity measure of the projected Poisson point process.

\begin{proposition}\label{prop:projetildef}
Let $\mu$ be a locally finite Borel measure on $(0,\infty)$, $\lambda>0$ and let
$\eta_{\mu,\lambda}$ be a Poisson point process on $\mathbb{H}^d\times(0,\infty)$
with intensity measure $\Lambda_{\mu,\lambda}$,
whose Laguerre diagram is almost surely a tessellation.
Then $\tau_k(\eta_{\mu,\lambda})$ is a Poisson point process on
$\HH^k\times(0,\infty)$ with intensity measure
$\Lambda_{\widetilde\mu,\lambda}$, where
$$
 \widetilde\mu(\mathrm dh)
  =\left[
     \omega_{d-k}\,h^k
     \int_{(0,h)}(h^2-t^2)^{\frac{d-k}{2}-1}
                   t^{1-d}\,\mu(\mathrm dt)
   \right]\mathrm dh.
$$
\end{proposition}
Lastly, we will show that if the weighted points arise from a Poisson point
process, then the section of the Laguerre diagram coincides almost
surely with the Laguerre diagram generated by the projected points.
For the Poisson processes considered here, we show that almost surely
every non-empty cell section has non-empty interior relative to $L_k$.
\begin{lemma}\label{lem:project2}
Let $\eta$ be a Poisson
point process with intensity measure either
$\Lambda_{\mu,\lambda}$ on $\mathbb{H}^d\times(0,\infty)$ or
$\Vol_{\SPH^{d-1}}\otimes(h^{d-2}\,\dd h)$ on
$\partial_{\infty}\mathbb{H}^d\times(0,\infty)$.
In the first case, let $\mu$ be a locally finite Borel
measure on $(0,\infty)$ and assume that the associated Laguerre diagram is almost
surely a tessellation. Then, almost surely,
$$
\mathcal L_{\rm hyp}^d (\eta)\cap L_k
=
\mathcal L_{\rm hyp}^k (\tau_k(\eta)).
$$
\end{lemma}

\begin{proof}[Proof of \cref{thrm:intersecdistrib}]
Since the distribution of the Poisson-Laguerre tessellation is invariant under isometries of $\mathbb{H}^d$, without loss of generality, we may assume that $L_k=\mathbb{H}^d\cap\tilde L_k,$ where $\tilde L_k=\{x\in\RR^{d+1}:x_{k+1}=\cdots=x_d=0\}$. 

Let $\xi_{\mu,\lambda}:=\tau_k(\eta_{\mu,\lambda})$. By Lemma~\ref{lem:project2} we have almost surely
$$
\mathcal{L}_{\rm hyp}^d(\eta_{\mu,\lambda})\cap L_k = \mathcal{L}_{\rm hyp}^{k}(\xi_{\mu,\lambda}).
$$  
By Proposition~\ref{prop:projetildef}, the process
$\xi_{\mu,\lambda}$ is a Poisson point process
with intensity measure
$\Lambda_{\widetilde\mu,\lambda}$.
Thus $\xi_{\mu,\lambda}\stackrel{d}{=}\eta_{\widetilde\mu,\lambda}$,
and the preceding identity yields
$$
\mathcal L_{\rm hyp}^{d}(\mu,\lambda)\cap L_k
\stackrel{d}{=}
\mathcal L_{\rm hyp}^{k}(\widetilde\mu,\lambda).
$$

The section inherits local finiteness and coverage from the original
tessellation, while Lemma~\ref{lem:project2} ensures that every nonempty
cell section has nonempty relative interior. Hence the resulting
diagram is almost surely a tessellation of $L_k$.
\end{proof}

The rest of the subsection is devoted to the proofs of the lemmata and propositions used.

\begin{proof}[Proof of Lemma~\ref{prop:project}] 
By definition of the polar polyhedron,
$$
P^d(\widehat{A})\cap \tilde{L}_k= \bigcap_{x\in \widehat{A}} (x^{\circ})^{+} \cap \tilde{L}_k.
$$
Fix $x\in \widehat{A}$ and $y\in\tilde L_k$. Since
$y_{k+1}=\cdots=y_d=0$, we have
\[
\loren{y}{x}
=
\sum_{i=1}^k y_ix_i-y_{d+1}x_{d+1}
=
\loren{y}{\proj_{\tilde L_k}(x)}.
\]
Consequently,
$ y\in(x^\circ)^+\cap\tilde L_k$ if and only if
$\loren{y}{x}=\loren{y}{\proj_{\tilde L_k}(x)}\leq-1$, which is equivalent to $y\in
((\proj_{\tilde L_k}(x))^\circ)^+$,
where the last polar halfspace is considered in $\tilde L_k$.
Therefore,
$$
P^d(\widehat{A})\cap\tilde L_k
=
\bigcap_{x\in \widehat{A}}
\left(
\left(\proj_{\tilde L_k}(x)\right)^\circ
\right)^+
=
\bigcap_{z\in \widehat{B}}(z^\circ)^+
=
P^k(\widehat{B}).$$
\end{proof}

\begin{proof}[Proof of Lemma~\ref{lem:coordexpre}]
After identifying $\tilde L_k$ with $\RR^{k+1}$, we have $\proj_{\tilde L_k}(\widehat{\alpha}_d(v,h))
=
h(v_1,\ldots,v_k,v_{d+1})$,
and $v_1^2+\cdots+v_k^2-v_{d+1}^2
=
\|v\|_L^2-\|\barbelow{v}\|^2$.

If $v\in\mathbb{H}^d$, then $\Vert v\Vert_L^2=-1$ and the projected vector is in $\mathcal{C}_{\mathrm{time}}^k$, then
$$
H(\proj_{\tilde L_k}(\widehat{\alpha}_d(v,h)))
=
h\sqrt{1+\|\barbelow{v}\|^2}.
$$
By the definition \eqref{def:phi} of $\widehat{\alpha}^{-1}_k$ we obtain that $\tau_k(v,h)=(v^{(k)}(v),h^{(k)}(v,h))$, where 
$$
v^{(k)}(v)=\frac{1}{\sqrt{1+\|\barbelow{v}\|^2} } (v_1,\ldots,v_k,v_{d+1})   
\quad 
\text{ and }
\quad 
h^{(k)}(v,h) = h\sqrt{1+\|\barbelow{v}\|^2} .
$$
If $v\in\partial_{\infty}\mathbb{H}^d$ and $\|\barbelow{v}\|>0$, then $\Vert v\Vert_L^2=0$ and the projected vector is again in $\mathcal{C}_{\mathrm{time}}^k$, with
$$
H(\proj_{\tilde L_k}(\widehat{\alpha}_d(v,h)))
=
h\|\barbelow{v}\|.
$$
By the definition \eqref{def:phi} of $\widehat{\alpha}^{-1}_k$ we obtain that $\tau_k(v,h)=(v^{(k)}(v),h^{(k)}(v,h))$, where 
$$
v^{(k)}(v)=\frac{1}{\|\barbelow{v}\| } (v_1,\ldots,v_k,v_{d+1})   
\quad 
\text{ and }
\quad 
h^{(k)}(v,h) =  h\|\barbelow{v}\|  .
$$
Finally, if $v\in\partial_{\infty}\mathbb{H}^d$ and $\|\barbelow{v}\|=0$, then $\Vert v\Vert_L^2=0$ and the projected vector is in $\mathcal{C}_{\mathrm{light}}$, with
$$
H(\proj_{\tilde L_k}(\widehat{\alpha}_d(v,h)))
=
h v_{d+1} = h.
$$
By \eqref{def:phi} we conclude that $\tau_k(v,h)=(v^{(k)}(v),h^{(k)}(v,h))$, where 
$$
v^{(k)}(v)=(v_1,\ldots,v_k,v_{d+1})   
\quad 
\text{ and }
\quad 
h^{(k)}(v,h) =  h .
$$

\end{proof}

\begin{proof}[Proof of Proposition~\ref{prop:projetildef}]
By the mapping theorem \cite[Theorem~24.16]{Klenke}, it is enough to determine how the intensity measure $\Lambda_{\mu,\lambda}$ of $\eta_{\mu,\lambda}$ changes under the map $\tau_k$.

Let $B\subset\HH^k$ be measurable and let $[s_1,s_2]\subset(0,\infty)$. To determine the push-forward measure, it suffices to compute its
values on rectangles $B\times[s_1,s_2]$, where $B\subset\HH^k$
is Borel and $0<s_1<s_2$, since these rectangles generate the
Borel $\sigma$-algebra of $\HH^k\times(0,\infty)$. Set
$$ I := \Lambda_{\mu,\lambda}\left(\tau_k^{-1}\bigl(B\times[s_1,s_2]\bigr)\right) = \lambda\iint_{\tau_k^{-1}(B\times[s_1,s_2])} \mu(\dd h)\,\Vol_{\mathbb{H}^d}(\dd v). $$
Using indicator functions, we can write
$$
I = \lambda\int_{\mathbb{H}^d}\int_{(0,\infty)} \mathbf1\left\{ \frac{(v_1,\ldots,v_k,v_{d+1})}{\sqrt{1+\|\barbelow{v}\|^2}}
\in B \right\} \mathbf{1}\left\{ s_1\leq h\sqrt{1+\|\barbelow{v}\|^2}\leq s_2 \right\} \mu(\dd h)\,\Vol_{\mathbb{H}^d}(\dd v).
$$

We now separate the first $k$ spatial coordinates from the remaining $d-k$ coordinates. Let $x:=(v_1,\ldots,v_k)$ and $y:=\barbelow{v}=(v_{k+1},\ldots,v_d)$, so that
$v_{d+1}=\sqrt{1+\|x\|^2+\|y\|^2}$. Since
$\Vol_{\mathbb{H}^d}(\dd v)=(1+\|x\|^2+\|y\|^2)^{-1/2}\,\dd x\,\dd y$,
we obtain
$$
I
=
\lambda\int_{\RR^k}\int_{\RR^{d-k}}\int_{(0,\infty)}
\mathbf1\left\{
\frac{(x,\sqrt{1+\|x\|^2+\|y\|^2})}
{\sqrt{1+\|y\|^2}}
\in B
\right\}
\frac{\mathbf1\left\{
s_1\leq h\sqrt{1+\|y\|^2}\leq s_2
\right\}}{\sqrt{1+\|x\|^2+\|y\|^2}}\mu(\dd h)\,\dd y\,\dd x.
$$
Define $u:=(1+\|y\|^2)^{-1/2}x$. Then
$1+\|x\|^2+\|y\|^2=(1+\|y\|^2)(1+\|u\|^2)$ and
$\dd x=(1+\|y\|^2)^{k/2}\,\dd u$.
Substituting these identities gives
$$
I
=
\lambda\int_{\RR^k}
\frac{\mathbf1\{(u,\sqrt{1+\|u\|^2})\in B\}}
{\sqrt{1+\|u\|^2}}
\int_{\RR^{d-k}}\int_{(0,\infty)}
\mathbf1\left\{
s_1\leq h\sqrt{1+\|y\|^2}\leq s_2
\right\}
(1+\|y\|^2)^{\frac{k-1}{2}}
\,\mu(\dd h)\,\dd y\,\dd u.
$$

The integrals separate, and the first factor becomes
$$
\int_{\RR^k}
\frac{\mathbf1\{(u,\sqrt{1+\|u\|^2})\in B\}}
{\sqrt{1+\|u\|^2}}\,\dd u
=
\Vol_{\HH^k}(B),
$$
which is precisely the hyperbolic volume of $B$ in the hyperboloid model. Thus,
$$
I
=
\lambda\Vol_{\HH^k}(B)
\int_{\RR^{d-k}}\int_{(0,\infty)}
\mathbf1\left\{
s_1\leq h\sqrt{1+\|y\|^2}\leq s_2
\right\}
(1+\|y\|^2)^{\frac{k-1}{2}}
\,\mu(\dd h)\,\dd y.
$$
  
We now express the remaining integral as a measure on the weight coordinate. For measurable $A\subset(0,\infty)$, set
$$
\widetilde\mu(A)
:=
\int_{(0,\infty)}\int_{\RR^{d-k}}
\mathbf1\left\{h\sqrt{1+\|y \|^2}\in A\right\}
(1+\|y \|^2)^{\frac{k-1}{2}}
\,\dd y\,\mu(\dd h).
$$
Changing the order of integration, we obtain
$$
I
=
\lambda\Vol_{\HH^k}(B)\,
\widetilde\mu([s_1,s_2]).
$$
It follows that the push-forward intensity measure is
$\Lambda_{\widetilde\mu,\lambda}$.

 To compute $\widetilde\mu$, we use spherical coordinates
$y=r\theta$, where $\theta\in\SPH^{d-k-1}$, and denote the original weight by $t$. Hence
$$
\widetilde\mu(A)
=
\omega_{d-k}
\int_{(0,\infty)}\int_0^\infty
\mathbf1\left\{t\sqrt{1+r^2}\in A\right\}
r^{d-k-1}(1+r^2)^{\frac{k-1}{2}}
\,\dd r\,\mu(\dd t).
$$
Finally, for fixed $t>0$, set $h=t\sqrt{1+r^2}$. Then
$r=t^{-1}\sqrt{h^2-t^2}$ and
$\dd r=\frac{h}{t}(h^2-t^2)^{-1/2}\,\dd h$.
As $r$ increases from $0$ to $\infty$, the variable $h$ increases
from $t$ to $\infty$. Therefore,
\begin{align*}
\widetilde\mu(A)
&=
\omega_{d-k}
\int_{(0,\infty)}\int_t^\infty
\mathbf1\{h\in A\}
h^k(h^2-t^2)^{\frac{d-k}{2}-1}t^{1-d}
\,\dd h\,\mu(\dd t)\\
&=
\int_A
\left[
\omega_{d-k}h^k
\int_{(0,h)}
(h^2-t^2)^{\frac{d-k}{2}-1}t^{1-d}
\,\mu(\dd t)
\right]\dd h.
\end{align*}
\end{proof}

\begin{proof}[Proof of Lemma~\ref{lem:project2}]
Set $\widehat{\eta}=\widehat{\alpha}_d(\eta)$ and $\xi=\widehat{\alpha}_k(\tau_k(\eta))$. Recall that due to Proposition \ref{prop:polyhedral} the Laguerre diagram can equivalently be defined as
$$
\mathcal{L}^{d}_{\rm hyp}(\eta )
= \big\{\proj_{\mathbb{H}^d}(F\cap\mathcal{C}^d_{\mathrm{time}}):
F\in\mathcal{F}_d(P(\widehat{\eta})),
\ F\cap\mathcal{C}^d_{\mathrm{time}}\neq\emptyset\big\}.
$$
Our first observation is that, for every
$F\in\mathcal{F}_d(P^d(\widehat{\eta}))$ with $F\cap\mathcal{C}^d_{\mathrm{time}}\neq\emptyset$ it holds 
\begin{equation}\label{eq:Projection}
\proj_{\mathbb{H}^d}(F\cap\mathcal{C}^d_{\mathrm{time}})\cap L_k = \proj_{\mathbb{H}^d}(F\cap\mathcal{C}^d_{\mathrm{time}}\cap \tilde{L}_k)=\proj_{\mathbb{H}^k}(F\cap\mathcal{C}^k_{\mathrm{time}}).
\end{equation}
Indeed, suppose that
$y\in\proj_{\mathbb{H}^d}(F\cap\mathcal{C}^d_{\mathrm{time}})\cap L_k$. Then there exists $x\in F\cap\mathcal{C}^d_{\mathrm{time}}$ such that
$\proj_{\mathbb{H}^d}(x)=y$. Since
$x=H(x)y$ and $y\in L_k\subset\tilde L_k$, the linearity of
$\tilde L_k$ implies that $x\in\tilde L_k$. Hence
$x\in F\cap\mathcal{C}^d_{\mathrm{time}}\cap\tilde L_k$, and therefore
$y\in\proj_{\mathbb{H}^d}(F\cap\mathcal{C}^d_{\mathrm{time}}\cap\tilde L_k)$.
Conversely, if
$y\in\proj_{\mathbb{H}^d}(F\cap\mathcal{C}^d_{\mathrm{time}}\cap\tilde L_k)$, then
$y\in\proj_{\mathbb{H}^d}(F\cap\mathcal{C}^d_{\mathrm{time}})$. Moreover, the radial projection of every
point of $\tilde L_k\cap\mathcal{C}^d_{\mathrm{time}}$ belongs to $L_k$. Hence $y\in L_k$,
which proves the reverse inclusion. 

The second observation is that, if
$F'\in\mathcal{F}_k(P^k(\xi))$ with $F'\cap \mathcal{C}_{\mathrm{time}}^k\neq \emptyset $, then there exists
$F\in\mathcal{F}_d(P^d(\widehat{\eta}))$ such that $F'=F\cap\tilde L_k$.
Indeed, since $F'$ is a a facet of $P^k(\xi)$, there exists
$y\in\xi$ such that $F'=y^\circ\cap P^k(\xi)$.
Choose $x\in\widehat{\eta}$ such that $y=\proj_{\tilde L_k}(x)$.
By Lemma~\ref{prop:project}, $P^k(\xi)=P^d(\widehat{\eta})\cap\tilde L_k$.
Moreover, for points in $\tilde L_k$, the polar hyperplane of
$y=\proj_{\tilde L_k}(x)$ is the intersection of $x^\circ$ with
$\tilde L_k$. Therefore,
$$
F' = y^\circ\cap P^k(\xi)= x^\circ\cap\tilde L_k\cap P^d(\widehat{\eta}) =
\left(x^\circ\cap P^d(\widehat{\eta})\right)\cap\tilde L_k.
$$
In particular, $F\cap\mathcal C^d_{\mathrm{time}}\neq\emptyset$,
so the corresponding Laguerre cell is nonempty.
This cell is full-dimensional by
Lemma~\ref{lem:full-dimensional-cells} in the interior case,
and by \cite[Theorem~3.2]{DAchilleEtAl2026} in the boundary case.
Hence $F\in\mathcal F_d(P^d(\widehat\eta))$  and
$F'=F\cap\tilde L_k$.

The final step uses general position of the Poisson process. We show that, almost surely, for every
$F\in \mathcal{F}_d(P^d(\widehat{\eta}))$ sattisfying $F\cap \tilde{L}_k\cap\mathcal{C}^d_{\mathrm{time} } = F\cap \mathcal{C}^k_{\mathrm{time}} \neq \emptyset$, the intersection
$F\cap\tilde L_k$ is a $k$-facet of $P^k(\xi)$. 

Suppose that $F\cap\tilde{L}_k$ is an $l$-dimensional face for some $l<k$, and set $S:=F\cap\tilde{L}_k$. Let $x_0\in\widehat{\eta}$ be the point generating the facet $F$ and let $y_0:=\proj_{\tilde{L}_k}(x_0)\in\xi$. Then $S=y_0^\circ\cap P^k(\xi)$. 
Choose $z\in S\cap\mathcal C^k_{\mathrm{time}}$.
By condition~\textup{(P1)}, in a neighbourhood of $z$,
$P^k(\xi)$ agrees with a $(k+1)$-dimensional polyhedron $Q$
defined by finitely many polar inequalities from $\xi$,
including the inequality generated by $y_0$.
The face $Q\cap y_0^\circ$ has dimension $l$ and is therefore
contained in at least $k+1-l$ facets of $Q$.
Let $y_1,\ldots,y_{k+1-l}\in\xi$ be points generating $k+1-l$ such facets. Define the linear subspace $\tilde I:=\operatorname{span}(S)\subset\tilde L_k$. Since $S$ has affine dimension $l$ and is contained in a polar hyperplane, which does not contain the origin, we have $\dim\tilde I=l+1$. For every $z\in S$ and every $j=1,\ldots,k+1-l$, we have
$\loren{z}{y_0}=-1$ and $\loren{z}{y_j}=-1$.
Therefore, $\loren{z}{y_j-y_0}=0$ for every $z\in S$. By linearity, this equality holds for every $z\in\tilde I$, and hence
$y_j-y_0\in\tilde I^\perp$, for $j=1,\ldots,k+1-l$.

Let $\tilde I^\perp$ denote the Lorentzian orthogonal complement of
$\tilde I$ inside $\tilde L_k$.  
However, $\dim\tilde I^\perp=k-l$, so the $k+1-l$ vectors $y_1-y_0,\ldots,y_{k+1-l}-y_0$ are linearly dependent. Equivalently, the projected points $y_0,y_1,\ldots,y_{k+1-l}$ are affinely dependent. 

The intensity measure of $\xi$ is absolutely continuous with respect
to Lebesgue measure on $\tilde L_k$. In the interior case, this follows
from Proposition~\ref{prop:projetildef}, since $\widetilde\mu$ has a
density and $\widehat{\alpha}$ is a smooth change of coordinates. In the boundary
case, the projected vectors have the form
$h(\theta_1,\ldots,\theta_k,1)$, where the projection of spherical
volume onto the first $k$ coordinates is absolutely continuous and
the weight measure has density $h^{d-2}$.
The multivariate Mecke formula \cite[Chapter~4]{LastPenrose}
therefore implies that almost surely no collection of at most $k+2$
distinct projected points is affinely dependent.
Consequently, almost surely, $F\cap\tilde L_k$ cannot have dimension smaller than $k$. Thus, whenever $F\cap\tilde L_k\cap\mathcal C^d_{\mathrm{time}}$ is non-empty, the intersection $F\cap\tilde L_k$ is a $k$-dimensional facet of $P^k(\xi)$.

Combining the preceding observations and omitting empty
intersections, by \eqref{eq:Projection} and Proposition \ref{prop:polyhedral} we obtain
\begin{align*}
\mathcal L_{\rm hyp}^{d}(\eta)\cap L_k
&=
\left\{
\proj_{\mathbb{H}^d}(F\cap\mathcal{C}^d_{\mathrm{time}}\cap\tilde L_k):
F\in\mathcal F_d(P^d(\widehat{\eta})),\
F\cap\mathcal{C}^d_{\mathrm{time}}\cap\tilde L_k\neq\emptyset
\right\}\\
&=
\left\{
\proj_{\HH^k}(F'\cap\mathcal{C}^k_{\mathrm{time}}):
F'\in\mathcal F_k(P^k(\xi)), F'\cap\mathcal{C}^k_{\mathrm{time}}\neq\emptyset
\right\}\\
&=
\mathcal L_{\rm hyp}^{k}(\tau_k(\eta)).
\end{align*}
Here we used that radial projection onto $\mathbb{H}^d$, restricted
to $\tilde L_k$, agrees with radial projection onto $\HH^k$
under the identification $\tilde L_k\cong\RR^{k+1}$.
\end{proof}

\subsection{Sections of Poisson-Voronoi tessellations}

Since Poisson-Voronoi tessellations are Poisson-Laguerre tessellations
with weight measure $\mu=\delta_1$, Theorem~\ref{thrm:intersecdistrib}
identifies their $k$-dimensional totally geodesic sections as
Poisson-Laguerre tessellations.

\begin{proof}[Proof of Theorem~\ref{thrm:Voronoiintersection}]
Applying Theorem~\ref{thrm:intersecdistrib} with $\mu=\delta_1$ gives
$$
\widetilde\mu(\dd h)
=
\omega_{d-k}h^k
(h^2-1)^{\frac{d-k}{2}-1}
\mathbf1\{h>1\}\,\dd h.
$$
Since $c_{\frac{d-k}{2}-1}=\omega_{d-k}$, this is precisely
$f_{k,\frac{d-k}{2}-1}(h)\,\dd h$. The intensity parameter is
preserved, and the proof follows.
\end{proof}

The following corollaries  follows directly from \cref{thrm:intersecdistrib}. They are hyperbolic analogues of \cite[Theorem 4.1 (i) and (iii)]{GusakovaKabluchkoThaele2024}.

\begin{corollary}\label{cor:beta} 
Let $d\geq2$, $\beta>-1$, and $1\leq k\leq d-1$. Let $\eta_{f_{d,\beta}}$ be the Poisson point process corresponding to the hyperbolic $\beta$-Voronoi intensity density $f_{d,\beta}$ defined by \eqref{eq:BetaDensity}. If $L_k$ is a $k$-dimensional subspace of $\mathbb{H}^d$, then 
$$\mathcal{L}^{d}_{\rm hyp}({f_{d,\beta}})\cap L_k \stackrel{d}{=} \mathcal{L}_{\rm hyp}^{k} \left( {f_{k,\beta+\frac{d-k}{2}}} \right). 
$$
\end{corollary}
\begin{proof} 
Since $f_{d,\beta}(h)=0$ for $h\leq 1$, the projected density
$\tilde f$ also vanishes on $(0,1]$. Let $h>1$. By
\cref{thrm:intersecdistrib},
$$
\tilde{f}(h)=\omega_{d-k} h^k \int_1^h  c_{\beta} (t^2-1)^\beta  
 (h^2-t^2)^{\frac{d-k}{2}-1} t   \, \dd t.
$$
Substituting $s=(t^2-1)/(h^2-1)$, with
$t\,\dd t=\frac12(h^2-1)\,\dd s$, yields
$$
\tilde f(h)
= \frac12\,\omega_{d-k}\,c_{\beta}\,h^{k}\,
(h^{2}-1)^{\,\beta+\frac{d-k}{2}}\,
\int_{0}^{1} s^{\beta}(1-s)^{\frac{d-k}{2}-1}\,ds.
$$
The integral on the right hand side is the Beta function $B(\beta+1,\frac{d-k}{2})$, and note that $B(a,b)=\frac{\Gamma(a)\Gamma(b)}{\Gamma(a+b)}$ and that $\omega_m = 2\pi^{m/2}/\Gamma \left(\frac{m}{2}\right)
$. So we obtain
$$
\frac12\omega_{d-k}c_\beta B\left(\beta+1,\frac{d-k}{2}\right) = \frac12\omega_{d-k}c_\beta \frac{\Gamma(\beta+1)\Gamma(\frac{d-k}{2})} {\Gamma(\beta+\frac{d-k}{2}+1)} = c_{\beta+\frac{d-k}{2}}.
$$
Thus $\tilde{f}(h)= f_{ k,\,\beta+\frac{d-k}{2} }(h)$ and the corollary follows.
\end{proof}

\begin{corollary}\label{cor:gaussian} 
Let $d\geq2$, $a>d-2$, and $1\leq k\leq d-1$. Let $\eta_{g_a}$ be the Poisson point process corresponding to the hyperbolic Gaussian-Voronoi intensity density $g_a(h)=h^a$. If $L_k$ is a $k$-dimensional subspace of $\mathbb{H}^d$, then 
$$\mathcal{L}^{d}_{\rm hyp}(g_a)\cap L_k \stackrel{d}{=} \mathcal{L}_{\rm hyp}^{k} (g_a). 
$$
\end{corollary}

\begin{proof}
    We note that
    \[
\tilde{g}(h)=\omega_{d-k} h^k \int_0^h   
 (h^2-t^2)^{\frac{d-k}{2}-1} t^{1-d+a}   \, \dd t=
\frac{\pi^{\frac{d-k}{2}}\Gamma(\frac{a-d+2}{2})}
{\Gamma(\frac{a-k+2}{2})}h^a,
\]
and the proof follows by combination of Theorem \ref{thrm:intersecdistrib} and Theorem \ref{thm:Gaussian}.
\end{proof}

\subsection{Sections of ideal Poisson-Voronoi tessellations}

We identify $\partial_\infty\HH^d$ with $\SPH^{d-1}$ and represent
$\widehat{\theta}\in\partial_\infty\HH^d$ by $(\theta,1)\in\partial_{\infty}\mathbb{H}^d$. Let $d\geq2$ and $1\leq k\leq d-1$. Let $\zeta$ be a Poisson
point process on $\partial_{\infty}\mathbb{H}^d\times(0,\infty)$ with
intensity measure
$$
\Lambda_\infty
=
\Vol_{\SPH^{d-1}}\otimes \, h^{d-2}\,\dd h.
$$

Before proving Theorem~\ref{thrm:idealsection}, we state a lemma giving the value of the integral needed in the proof. We prove the lemma after the theorem.
\begin{lemma}\label{le:calculationIntegral}
Let $d\geq 2$ and $1\leq k \leq d-1$. Let $B\subset \HH^k\times (0,\infty)$ be a Borel set, then
\[
\iint_{\tau_k^{-1}(B)}  h^{d-2} \, \dd h \, \Vol_{\SPH^{d-1}} (\dd \theta)= \omega_{d-k} \iint_B h^{d-2} \, \Vol_{\HH^k}(\dd v) \, \dd h.
\]
\end{lemma}

\begin{proof}[Proof of Theorem~\ref{thrm:idealsection}]
We begin by recalling that the distribution of $\mathcal{V}_{\infty}$ is invariant under isometries (see \cite{DAchilleEtAl2026}). 
 Hence, without loss of generality, it suffices to consider $L_k=\HH^d\cap  \tilde{L}_k$ 
with $\tilde{L}_k=\{(x_1,\ldots,x_{d+1}) : x_{k+1}=\ldots=x_d=0\}$. 
By Lemma~\ref{lem:project2}, almost surely, 
$$
\mathcal V_\infty \cap L_k= \mathcal L_{\rm hyp}^d (\zeta)\cap L_k
=
\mathcal L_{\rm hyp}^k (\xi),
$$
where 
$\xi:=\tau_k(\zeta)$. Note, that since $\partial_\infty L_k\cong\SPH^{k-1}$ has zero spherical
volume in $\partial_\infty\HH^d$, almost surely no nucleus of $\zeta$ lies
in $\partial_\infty L_k$. Thus every point of $\proj_{\tilde L_k}(\widehat{\alpha}(\zeta))$ is in $\mathcal{C}^k_{\mathrm{time}}$,
and its location and weight coordinates belong to
$\HH^k\times(0,\infty)$.

By the mapping theorem \cite[Theorem~24.16]{Klenke}, the image $\xi$ is a Poisson point process and its intensity measure is given by the push-forward 
\[
\Lambda_k(B)=\iint_{\tau_k^{-1}(B)} h^{d-2} \, \dd h \, \Vol_{\SPH^{d-1}} (\dd \theta).
\]
By Lemma~\ref{le:calculationIntegral}, the weight measure of
$\xi$ is $\omega_{d-k}h^{d-2}\,\dd h$. Since this is the intensity measure of a hyperbolic Gaussian tessellation, multiplying all weights
by $\omega_{d-k}^{1/(d-1)}$ removes the factor $\omega_{d-k}$
from the intensity and leaves the cells unchanged as follows by Theorem \ref{thm:Gaussian}. This proves
the claimed distributional identity. 
\end{proof}

 \begin{proof}[Proof of Lemma~\ref{le:calculationIntegral}]
Since $(\theta_1,\ldots,\theta_d)\in\mathbb{S}^{d-1}$, we have $\sum_{i=1}^{k} \theta_i^2 + \Vert \barbelow{\theta}\Vert^2 =1$, where $\barbelow{\theta}=(\theta_{k+1},\ldots,\theta_d)$.
The set $\{\theta\in\SPH^{d-1}:\barbelow{\theta}=0\}$ has zero spherical volume, since it is a $(k-1)$-dimensional
subsphere and $k<d$. We may therefore restrict the integral to
$\|\barbelow{\theta}\|>0$. By Lemma~\ref{lem:coordexpre} on this set we have
$$
\tau_k(\widehat{\theta},h)
=
\left(
\frac{(\theta_1,\ldots,\theta_k,1)}
     {\|\barbelow{\theta}\|},
h \|\barbelow{\theta}\|
\right).
$$
Introduce a change of variables
$$
v_i:=\frac{\theta_i}{\Vert \barbelow{\theta}\Vert},\,  \text{ for } i=1,\ldots,k, \quad v_{k+1}:=\frac{1}{\Vert \barbelow{\theta}\Vert }, \quad \text{ and } \quad \omega:=\left( \frac{\theta_{k+1}}{\Vert \barbelow{\theta} \Vert },\ldots, \frac{\theta_d}{\Vert\barbelow{\theta}\Vert}   \right).
$$
We also write $v:=(v_1,\ldots,v_{k+1})$ and $\bar{v} := (v_1,\ldots,v_k)$.
Dividing
$\sum_{i=1}^k\theta_i^2+\|\barbelow{\theta}\|^2=1$ by
$\|\barbelow{\theta}\|^2$ yields $\Vert \barbelow{\theta} \Vert = (1+\Vert \bar{v}\Vert^2)^{-1/2}$.
Consequently, we may rewrite $\theta=\left( \Vert \barbelow{\theta} \Vert \cdot \bar v,  \Vert \barbelow{\theta} \Vert  \cdot\omega \right) =  (1+\Vert \bar{v}\Vert^2 )^{-1/2} (\bar v,\omega)$,
where $v\in \HH^k$ and $\omega\in\mathbb{S}^{d-k-1}$.

First we decompose $\theta \in \mathbb{S}^{d-1}$ splitting the coordinates, $
\theta= (\sin\psi \, \beta, \cos \psi \, \omega)$,
where $\beta \in \mathbb{S}^{k-1}$, and $\omega\in \mathbb{S}^{d-k-1}$.
The surface measure then admits the decomposition
\begin{equation}\label{eq:changecoord2}
\Vol_{\SPH^{d-1}}(\dd \theta) = (\sin\psi)^{k-1} (\cos \psi )^{d-k-1} \, \dd \psi \, \Vol_{\SPH^{k-1}}(\dd \beta ) \Vol_{\SPH^{d-k-1}}(\dd \omega).
\end{equation}
Next, observe that for $i=1,\ldots,k$, we can write $v_i=\frac{\theta_i}{\Vert \barbelow{\theta}\Vert}= \frac{\sin\psi}{\cos\psi}\beta_i= (\tan\psi ) \beta_i$.
Since $\bar v=\rho\beta$ with $\rho=\tan\psi$, polar coordinates
in $\RR^k$ give
$$
\dd v_1\cdots\dd v_k
=
\rho^{k-1}\,\dd\rho\,\Vol_{\SPH^{k-1}}(\dd\beta)
=
(\tan\psi)^{k-1}(\sec\psi)^2\,\dd\psi\,
\Vol_{\SPH^{k-1}}(\dd\beta).
$$
Substituting into \eqref{eq:changecoord2} and using
$\cos\psi=(1+\|\bar v\|^2)^{-1/2}$, we obtain
\begin{align*}
\Vol_{\SPH^{d-1}}(\dd\theta)
&=
(\cos\psi)^d\,\dd v_1\cdots\dd v_k\,
\Vol_{\SPH^{d-k-1}}(\dd\omega)\\
&=
(1+\|\bar v\|^2)^{-d/2}\,\dd v_1\cdots\dd v_k\,
\Vol_{\SPH^{d-k-1}}(\dd\omega).
\end{align*}

For fixed $\bar v$, we change the weight coordinate by setting $t=h\|\barbelow{\theta}\|
=(1+\Vert \bar{v}\Vert^2 )^{-1/2} h$.
Thus $h=t(1+\Vert \bar{v}\Vert^2 )^{1/2}$ and $h^{d-2}\,\dd h = t^{d-2}(1+\|\bar v\|^2)^{(d-1)/2}\,\dd t$.
Combining the spherical volume decomposition with the change of
weight coordinate, we obtain
$$
\iint_{\tau_k^{-1}(B)}
h^{d-2}\,\dd h\,\Vol_{\SPH^{d-1}}(\dd\theta)
=
\int_{\tilde B}
t^{d-2}(1+\|\bar v\|^2)^{-1/2}
\,\dd t\,\dd v_1\cdots\dd v_k\,
\Vol_{\SPH^{d-k-1}}(\dd\omega),
$$
where $\tilde B: = \{(\bar v,\omega,t):(v,t)\in B\}$, and  $v=(\bar v,\sqrt{1+\|\bar v\|^2})$.
 
Since the condition $(v,h)\in B$ is independent of $\omega$,
integration over $\SPH^{d-k-1}$ gives the factor $\omega_{d-k}$.
Moreover, $v_{k+1}^2=1+\|\bar v\|^2$, so
$$
(1+\|\bar v\|^2)^{-1/2}\,\dd v_1\cdots\dd v_k
=
\Vol_{\HH^k}(\dd v).
$$
Consequently,
$$
\iint_{\tau_k^{-1}(B)}
h^{d-2}\,\dd h\,\Vol_{\SPH^{d-1}}(\dd\theta)
=
\omega_{d-k}\iint_B h^{d-2}\,\dd h\,\Vol_{\HH^k}(\dd v),
$$
as claimed.
\end{proof}

\addtocontents{toc}{\protect\setcounter{tocdepth}{0}}
\section*{Acknowledgments}
  This work was funded by the Deutsche Forschungsgemeinschaft (DFG, German Research Foundation) under Germany’s Excellence Strategy (EXC 2044/2 – 390685587, Mathematics Münster: Dynamics–Geometry–Structure) and through SPP 2265 Random Geometric Systems. P. S. acknowledges financial support from the FIS 2 project “Understanding pattern formation in nature via complex analysis” CUP B53C24009510001. 

  \section*{Declarations}

\subsection*{Statement on the use of generative AI}
The project was initiated in 2025, and the authors finalized most of the mathematical statements and proofs by March 2026 without using generative AI models. Theorem \ref{thrm:idealsection} and related results were proved by May 2026, also without the use of generative AI models. All theorem statements and proofs were developed by the authors. During the final stage of preparing the manuscript and bringing the results together, ChatGPT was used to improve the clarity and wording of the text, make it more concise, identify potential mathematical issues, and correct typographical errors. The authors take full responsibility for the contents of the paper.

\subsection*{Conflict of interest statement}
The authors declare that they have no conflict of interest.

\subsection*{Data availability statement}
We do not analyze or generate any datasets.

\addtocontents{toc}{\protect\setcounter{tocdepth}{1}}

\footnotesize{
\bibliographystyle{abbrv}
\bibliography{biblio}
}

\end{document}